\documentclass[11pt]{article}
\usepackage{amsmath,amssymb,amsthm}
\usepackage{mathtools}
\usepackage{url}
\usepackage{bbm}
\newcommand{\one}{\mathbbm{1}}
\usepackage{geometry}
\usepackage{tikz}
\usetikzlibrary{arrows.meta,positioning,calc}
\usepackage{caption}
\usepackage{subcaption}
\usepackage{xcolor}
\usepackage[normalem]{ulem}

\makeatletter
\newcommand{\xrightleftarrows}[2]{
  \mathrel{\mathop{
    \vcenter{\offinterlineskip\m@th
      \ialign{\hfil##\hfil\cr
        \hphantom{$\scriptstyle\mspace{8mu}{#1}\mspace{8mu}$}\cr
        \rightarrowfill\cr
        \vrule height0pt width 2em\cr
        \leftarrowfill\cr
        \hphantom{$\scriptstyle\mspace{8mu}{#2}\mspace{8mu}$}\cr
        \noalign{\kern-0.3ex}
      }
    }
  }\limits^{#1}_{#2}}
}
\makeatother

\newtheorem{theorem}{Theorem}[section]
\newtheorem{lemma}[theorem]{Lemma}
\newtheorem{proposition}[theorem]{Proposition}
\newtheorem{corollary}[theorem]{Corollary}
\newtheorem{conjecture}[theorem]{Conjecture}
\theoremstyle{definition}
\newtheorem{definition}[theorem]{Definition}
\newtheorem{example}[theorem]{Example}
\theoremstyle{remark}
\newtheorem{remark}[theorem]{Remark}
\newcommand{\exampleend}{\unskip\nobreak\hfill\ensuremath{\triangle}}

\title{Substrate hypergraph cycles:\\
unique steady states and persistence for every enzyme assignment;\\
oscillations enabled by enzyme sharing}
\author{Badal Joshi\thanks{California State University San Marcos, San Marcos,
California, USA.}
\and Tung D. Nguyen\thanks{Fulbright University Vietnam, Ho Chi Minh City,
Vietnam.}}
\date{}

\begin{document}
\maketitle

\begin{abstract}
We associate a directed substrate hypergraph with an enzymatic reaction
system.  The hypergraph records which substrates each composite reaction consumes and
produces, and each composite reaction is realized by a Henri--Michaelis--Menten (HMM)
module whose enzyme may be distinct, shared with other reactions, or absent.
We call the hypergraph a cycle when it is strongly connected and every
substrate is produced by exactly one reaction and consumed by exactly one
reaction, each time with coefficient one; equivalently, its dual is a strongly
connected ordinary directed graph.  For a cycle, every HMM realization, under
every enzyme assignment and for all positive rate constants, has exactly one positive
steady state in each stoichiometric compatibility class that contains a
strictly positive point, and no steady state on the boundary of such a class, and every trajectory in such a class is persistent.  The interior dynamics, in contrast, depend on how enzymes are
shared.  When
each enzyme catalyzes at most one reaction, every trajectory in such a class
converges to its unique steady state.  With enzyme sharing, convergence can
fail even for the three-substrate transformation cycle
$X_1\to X_2\to X_3\to X_1$: if two transformations share an enzyme and the third has
its own enzyme or none, the steady state can undergo a supercritical Hopf bifurcation that produces a stable limit cycle. Thus the substrate structure controls existence,
uniqueness, and persistence, while the pattern of enzyme sharing controls the
interior dynamics.  The results extend to weighted cycles, which allow higher stoichiometric coefficients,
and to longer chains of intermediates.
\end{abstract}



\section{Introduction}\label{sec:intro}

Many biochemical networks are known through their substrate transformations,
while the enzymes that catalyze them, the sharing of enzymes between reactions,
and the intermediate mechanisms are uncertain~\cite{SalazarHofer2009,Gunawardena2014}.
In human cells, for example, mass spectrometry had identified more than
100{,}000 phosphorylation sites by 2019, but kinases had been assigned to fewer
than 5\% of them, and a site may be targeted by a single kinase or by several
kinases in different contexts~\cite{NeedhamEtAl2019}; predicting the
responsible kinase from substrate data is itself an active
problem~\cite{JohnsonEtAl2023}.
A kinase may act on several
sites, a phosphatase on several phosphoforms, and a catalytic step may pass
through one or several intermediate complexes.  Adding intermediates or
sharing enzymes can preserve some properties of a network and create or
destroy others~\cite{FeliuWiuf2012,FeliuWiuf2013,SuwanmajoKrishnan2015,SuwanmajoRameshKrishnan2020};
for example, Feliu and Wiuf showed that two one-site futile cycles coupled only
through a shared kinase and a shared phosphatase can be multistationary,
although each cycle alone has a unique steady state in every
class~\cite{FeliuWiuf2012} (Example~\ref{ex:feliu_wiuf}).  This raises two questions:
which dynamical properties are determined by the substrate transformations
alone, for every enzymatic realization, and which depend in addition on how
enzymes are assigned and shared.  For a class of substrate hypergraphs that we
call cycles, we show that the substrate structure determines existence and
uniqueness of positive steady states as well as persistence, while convergence to
the steady state depends in addition on whether enzymes are shared.  The
substrate-hypergraph framework~\cite{JoshiNguyen2026Hypergraph} records the
substrate transformations and suppresses catalytic detail.

\begin{figure}[h!]
\centering
\begin{tikzpicture}[baseline=(current bounding box.center)]
  \node[align=center,font=\small] (mech) at (0,0) {$E+X\;\rightleftarrows\;C_1\;\longrightarrow\;E+Y$\\[2pt]
                                              $F+Y\;\rightleftarrows\;C_2\;\longrightarrow\;F+X$};
  \draw[-{Stealth[length=10pt,width=7pt]},line width=1.25pt] (3.1,0) -- (4.65,0)
        node[midway,above=13pt,font=\scriptsize]{erase enzymes and intermediates};
  \begin{scope}[shift={(5.55,0)}]
    \node[circle,draw,thick,minimum size=0.72cm,inner sep=1pt] (X) at (0,0) {$X$};
    \node[circle,draw,thick,minimum size=0.72cm,inner sep=1pt] (Y) at (2.0,0) {$Y$};
    \draw[-{Stealth[length=8pt,width=5.5pt]},line width=1.35pt] (X) to[bend left=24]
          node[above,font=\scriptsize\itshape]{$\rho_1$} (Y);
    \draw[-{Stealth[length=8pt,width=5.5pt]},line width=1.35pt] (Y) to[bend left=24]
          node[below,font=\scriptsize\itshape]{$\rho_2$} (X);
  \end{scope}
\end{tikzpicture}
\caption{From enzymatic mechanism to substrate skeleton.}
\label{fig:intro_skeleton}
\end{figure}
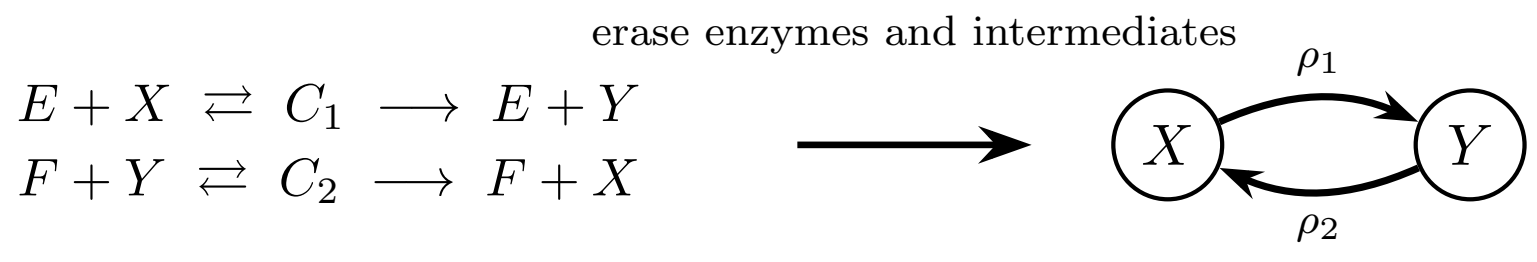

The one-site futile cycle is the elementary enzymatic system
\begin{equation}\label{eq:intro_futile}
  E+X \;\xrightleftarrows{k_1^+}{k_1^-}\; C_1 \;\xrightarrow{\,k_1\,}\; E+Y,
  \qquad
  F+Y \;\xrightleftarrows{k_2^+}{k_2^-}\; C_2 \;\xrightarrow{\,k_2\,}\; F+X .
\end{equation}
It has a unique positive steady state in every positive stoichiometric
compatibility class.  Erasing the enzymes and intermediates leaves the
\emph{substrate skeleton}, which records only which substrates each composite
reaction consumes and produces; for~\eqref{eq:intro_futile} it has two
edges (Figure~\ref{fig:intro_skeleton}).

We represent the skeleton as a directed hypergraph so that a
composite reaction may consume and produce several substrates, as in
phosphotransfer,
\begin{equation}\label{eq:intro_phosphotransfer}
  X_p+Y\longrightarrow X+Y_p,
\end{equation}
which cannot be represented by a single directed edge on substrate species.

\begin{figure}[h!]
\centering
\begin{subfigure}[b]{0.48\textwidth}
\centering
\begin{tikzpicture}[
  nd/.style={circle,draw,thick,minimum size=0.85cm,inner sep=1pt,font=\footnotesize},
  en/.style={rounded corners=3pt,draw,thick,minimum width=0.95cm,minimum height=0.5cm,inner sep=1.5pt,font=\scriptsize},
  AH/.style={-{Stealth[length=8pt,width=5.5pt]},line width=1.45pt},
  LN/.style={line width=1.45pt},
  lbl/.style={font=\scriptsize\itshape,inner sep=1pt},
  scale=0.82
]
  \node[nd] (X)   at (-2.1,0) {$X$};
  \node[nd] (Xp)  at ( 2.1,0) {$X_p$};
  \node[circle,draw,thick,fill=yellow!45,inner sep=0.6pt,font=\tiny]
    at ($(Xp.north east)+(0.02,0.02)$) {P};
  \node[en,fill=red!18]  (ATP) at (-1.4,2.5) {ATP};
  \node[en,fill=blue!12] (ADP) at ( 1.4,2.5) {ADP};
  \coordinate (J) at (0,1.2);
  \draw[LN,orange!85!red] (X) to[out=60,in=180] (J);
  \draw[AH,orange!85!red] (J) to[out=0,in=120] (Xp);
  \draw[LN,orange!85!red] (ATP) to[out=-90,in=180] (J);
  \draw[AH,orange!85!red] (J) to[out=0,in=-90] (ADP);
  \fill[orange!85!red] (J) circle (3pt);
  \node[lbl,text=orange!85!red,anchor=east] at (-1.45,1.45) {$\rho_1$, kinase};
  \draw[AH,teal!80!black] (Xp) to[out=-120,in=-60]
    node[lbl,below,text=teal!80!black] {$\rho_2$, phosphatase} (X);
  \draw[AH,violet!75!black,line width=2.2pt] (ADP) to[out=90,in=90,looseness=1.1]
    node[lbl,above,text=violet!75!black] {$\rho_3$, regeneration} (ATP);
\end{tikzpicture}
\caption{The substrate hypergraph.}
\end{subfigure}
\hfill
\begin{subfigure}[b]{0.46\textwidth}
\centering
\begin{tikzpicture}[
  rn/.style={circle,draw,line width=1.35pt,minimum size=0.86cm,inner sep=1pt,font=\small\itshape},
  el/.style={font=\scriptsize,fill=white,inner sep=1pt},
  AH/.style={-{Stealth[length=9pt,width=6pt]},line width=1.45pt},
  scale=1.1
]
  \node[rn,draw=orange!85!red,text=orange!85!red] (r1) at (0,0) {$\rho_1$};
  \node[rn,draw=teal!80!black,text=teal!80!black] (r2) at (-2.4,0) {$\rho_2$};
  \node[rn,draw=violet!75!black,text=violet!75!black] (r3) at (2.4,0) {$\rho_3$};
  \draw[AH] (r1) to[bend right=55] node[el] {$X_p$} (r2);
  \draw[AH] (r2) to[bend right=55] node[el] {$X$} (r1);
  \draw[AH] (r1) to[bend left=55] node[el] {ADP} (r3);
  \draw[AH] (r3) to[bend left=55] node[el] {ATP} (r1);
\end{tikzpicture}
\vspace{0.8em}
\caption{Its dual graph.}
\end{subfigure}
\caption{A kinase--phosphatase pair coupled to ATP regeneration has a substrate hypergraph which is a cycle
(Example~\ref{ex:atp_cycle}).  Each loop in~(a) is a directed cycle in~(b).
Colors distinguish reactions; branches joined at a dot are consumed or
produced together.}
\label{fig:intro_atp}
\end{figure}
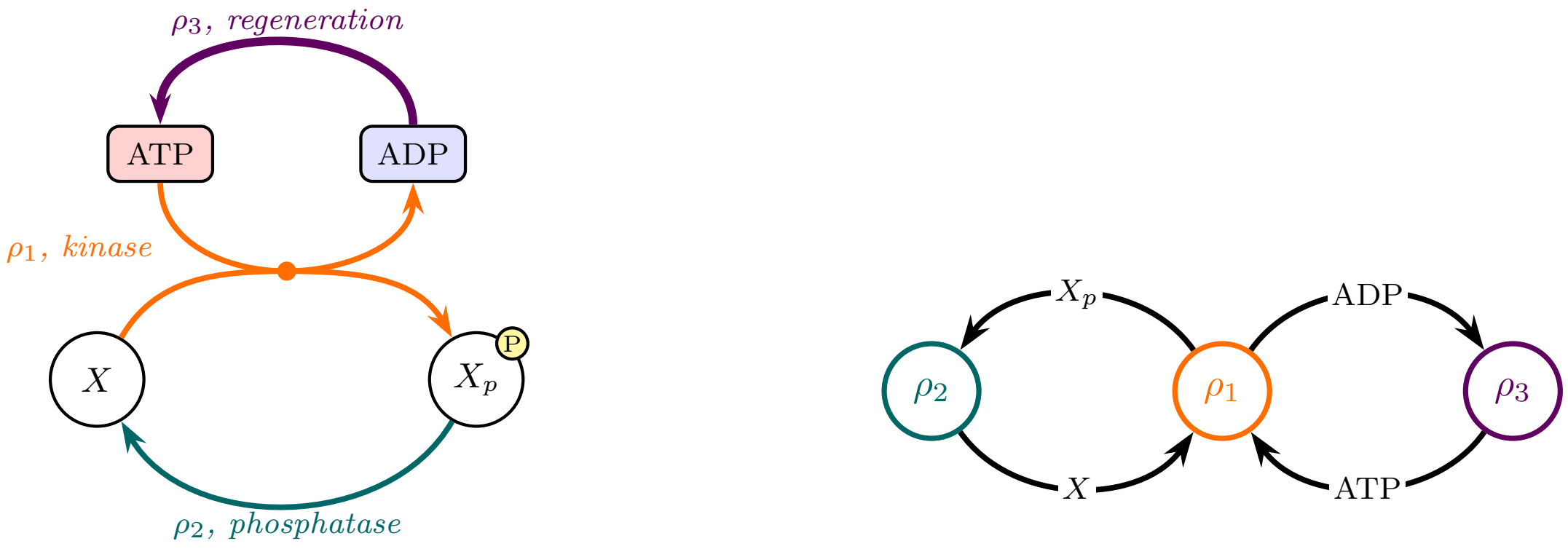

A substrate hypergraph is a \emph{cycle} if it is strongly connected and every
substrate occurs with total multiplicity one among both reaction tails and
reaction heads.  Equivalently, its dual is an ordinary strongly connected
directed graph whose edges are the substrates and whose nodes are the
reactions, each substrate pointing from the reaction that produces it to the
reaction that consumes it.  A ring of phosphotransfer reactions~\eqref{eq:intro_phosphotransfer}
is a cycle (Example~\ref{ex:transfer_cycle}), and so is a kinase--phosphatase
pair coupled to ATP regeneration (Example~\ref{ex:atp_cycle}; Figure~\ref{fig:intro_atp}).

An \emph{enzyme assignment} specifies, for each composite reaction
$\rho\colon t_\rho\to h_\rho$, which enzyme catalyzes it, if any; enzymes may
be distinct, shared between reactions, or absent.  In the \emph{HMM
realization} of an enzyme assignment, the reaction steps of each composite
reaction follow the one-intermediate Henri--Michaelis--Menten (HMM) module
\begin{equation}\label{eq:intro_module}
  E+t_\rho\;\rightleftarrows\;C_\rho\;\longrightarrow\;E+h_\rho,
\end{equation}
where $E$ is the assigned enzyme and is omitted for an enzyme-free reaction,
and the realization is governed by mass-action kinetics $\dot u=f(u)$, where
$u$ is the vector of species concentrations
(Section~\ref{sec:enzyme_assign}).  Unless stated otherwise, every result in
this paper concerns HMM realizations; Section~\ref{sec:extensions} also
treats chains of intermediates.  We call a stoichiometric compatibility class \emph{positive} if it contains a
strictly positive point.  With this terminology, our main result is the
following.
\begin{theorem}[Main theorem]\label{thm:marked_main}
Let $\mathcal H$ be a finite substrate hypergraph cycle.  For every enzyme
assignment and all positive rate constants, the HMM realization of
$\mathcal H$ has the following properties:
\begin{enumerate}
\item each positive stoichiometric compatibility class contains exactly one
  positive steady state and no steady state on the boundary of the
  nonnegative orthant;
\item every trajectory that starts in a positive stoichiometric
  compatibility class is strictly positive for all $t>0$, and its
  $\omega$-limit set lies in the open orthant; in particular, the
  realization is persistent;
\item $f(u^{(1)})\ne f(u^{(2)})$ for any two distinct strictly positive
  states $u^{(1)},u^{(2)}$ in the same stoichiometric compatibility class;
  that is, the realization is \emph{injective}.
\end{enumerate}
\end{theorem}

Theorem~\ref{thm:marked_main} deliberately makes no claim about stability or
convergence, and no such claim holds uniformly over enzyme assignments.  The
same substrate structure does yield a stronger dynamical conclusion when
enzymes are not shared.
\begin{theorem}[Global convergence without enzyme sharing]
\label{thm:distinct_convergence}
Let $\mathcal H$ be a finite substrate hypergraph cycle, and consider an
enzyme assignment in which each enzyme catalyzes at most one reaction;
enzyme-free reactions are allowed.  Then, for all positive rate constants,
every trajectory of the HMM realization that starts in a positive
stoichiometric compatibility class converges to the unique positive steady
state of that class.
\end{theorem}
  With enzyme sharing, convergence can fail.
\begin{theorem}[The smallest oscillating transformation cycle]
\label{thm:hopf_main}
The transformation cycles $X_1\to X_2\to X_3\to X_1$ with one enzyme
shared by two of their reactions and the remaining reaction either catalyzed
by a second enzyme or enzyme-free are minimal transformation cycles whose HMM realizations can have trajectories that
do not converge to the positive steady state:
\begin{enumerate}
\item for each of these two assignments, there exist positive rate constants and a positive compatibility class
  whose unique positive steady state is unstable;
\item for each of these two assignments, along a one-parameter family of rate constants, this steady state
  undergoes a supercritical Hopf bifurcation, and a limit cycle exists that is orbitally asymptotically stable within the compatibility class;
\item on a transformation cycle with at most three reactions, under any enzyme
  assignment with enzyme-free reactions allowed, some positive rate constants
  and some positive compatibility class admit a trajectory that does not
  converge to the positive steady state of that class if and only if the cycle
  has three reactions and some enzyme catalyzes exactly two of them.
\end{enumerate}
\end{theorem}

The dual graph drives all three results.  Its directed cycles carry
nonnegative conservation laws (Corollary~\ref{cor:dcycle_law}) that give uniqueness, boundary
exclusion, and the description of positive classes.  Its edges produce the
positive off-diagonal entries of the Jacobian in reaction coordinates, so that
with distinct enzymes the dynamics in these coordinates are cooperative and
irreducible (Figure~\ref{fig:influence}(a)), which is the hypothesis of the
convergence theorem of Angeli, De~Leenheer, and
Sontag~\cite{AngeliDeLeenheerSontag2010}.  When two reactions share an
enzyme, their binding steps compete for it and produce negative entries; the
network of Theorem~\ref{thm:hopf_main} has exactly this sign pattern
(Figure~\ref{fig:influence}(b)).  Thus the enzyme-incidence pattern, which records which reactions use
or share which enzymes, controls the interior dynamics.

\subsection{Related work}\label{sec:related}

\paragraph{Frameworks.}
Theorem~\ref{thm:marked_main} extends the substrate-hypergraph framework
of~\cite{JoshiNguyen2026Hypergraph} from existence to uniqueness and
persistence.  In that framework, currents constrain the existence of
positive steady states and bifunctional enzymes produce absolute
concentration robustness.  The post-translational modification framework of Thomson and
Gunawardena~\cite{ThomsonGunawardena2009,ThomsonGunawardena2009Nature} and the
MESSI framework of P\'erez Mill\'an and Dickenstein~\cite{PerezMillanDickenstein2018}
start from a given network and study it through a graph on the modification
forms of its substrates.  Transformation cycles in which every reaction is
catalyzed by an enzyme satisfy the hypotheses
of~\cite[Theorem~3.15]{PerezMillanDickenstein2018}, which gives persistence
and the absence of boundary steady states; Theorem~\ref{thm:marked_main} adds
uniqueness.  In these frameworks enzymatic reactions are one-to-one, the only
two-substrate reactions are the swaps of
MESSI~\cite[Definition~2.3]{PerezMillanDickenstein2018}, which
\cite[Remark~3.16]{PerezMillanDickenstein2018} excludes from Theorem~3.15,
donors such as ATP are absorbed into the rate
constants~\cite{ThomsonGunawardena2009Nature,PerezMillanDickenstein2018}, and
each conclusion concerns a fixed network rather than the class of its
enzymatic realizations.  By contrast, the hypergraph framework
of~\cite{JoshiNguyen2026Hypergraph} and of this paper allows binding,
splitting, and transfer reactions with arbitrary numbers of substrates and
arbitrary stoichiometric coefficients, such as $C\to A+B$ and $2A\to B+C$,
provided that each side of a reaction contains at least one substrate, so
flow reactions are excluded.  Currents were defined
in~\cite{JoshiNguyen2026Hypergraph}; current cones and hypergraph cycles are
defined here.

\paragraph{Steady states (Theorem~\ref{thm:marked_main}).}
For individual reaction networks, deficiency theory, injectivity and
concordance criteria, and the theory of minimal multistationary subnetworks
decide uniqueness or multistationarity, in many cases, from the detailed
reaction
graph~\cite{Feinberg1987,Feinberg1995,Feinberg2019,CraciunFeinberg2005,
MullerEtAl2016,BanajiPantea2016,ShinarFeinberg2012,JoshiShiu2012,
JoshiShiu2013,JoshiShiu2015,JoshiShiu2017}.  For two phosphorylation networks
whose skeletons are transformation cycles, uniqueness has been established
network by network: the cyclic multisite modification network of
Suwanmajo, Ramesh, and Krishnan~\cite{SuwanmajoRameshKrishnan2020} on four phosphoforms, for the common and separate
kinase and phosphatase variants they consider; and the mixed mechanism of Suwanmajo and
Krishnan~\cite{SuwanmajoKrishnan2015}, with processive phosphorylation and
distributive dephosphorylation~\cite{SuwanmajoKrishnan2015,ConradiMinchevaShiu2019}.
The uniqueness result for the mixed mechanism confirmed half of a conjecture
of Conradi and Shiu~\cite{ConradiShiu2015}.  Both are special cases of
Theorem~\ref{thm:marked_main} and its extension to chains of intermediates
(Corollary~\ref{cor:special_cases}).  Distributive multiple futile cycles, whose
hypergraph skeletons are not cycles, can be multistationary, as Wang and Sontag showed~\cite{WangSontag2008}, and the
number of their stable steady states can grow with the number of
sites~\cite{ThomsonGunawardena2009Nature}.

\paragraph{Convergence (Theorem~\ref{thm:distinct_convergence}).}
Global convergence was proved for the one-site futile cycle by Angeli and
Sontag~\cite{AngeliSontag2008}, for processive multisite phosphorylation by
Conradi and Shiu~\cite{ConradiShiu2015}, for a class of cyclic processive systems by
Rao~\cite{Rao2017}, and, for transformation cycles
$X_1\to X_2\to\cdots\to X_m\to X_1$ in which each reaction has its own enzyme
and a chain of intermediates, by Eithun and Shiu~\cite{EithunShiu2017}.  Their
proof, like ours in Section~\ref{sec:convergence}, combines the convergence
theorem of~\cite{AngeliDeLeenheerSontag2010} with the intermediate-removal
results of~\cite{MarcondesDeFreitasWiufFeliu2017Convergence}.
Theorem~\ref{thm:distinct_convergence} extends this from transformation
cycles to hypergraph cycles and allows enzyme-free reactions; the results
above are special cases (Corollary~\ref{cor:special_cases}).
Donnell and Banaji~\cite{DonnellBanaji2013} proved global stability under
general kinetics for networks whose stoichiometric matrix has a special
factorization; the enzymatic futile cycle is one of their examples.
Transformation cycles with distinct enzymes have such a factorization, but not, in general, when an enzyme is shared, when a reaction has two tail substrates, or
for chains of intermediates.

\paragraph{Oscillations (Theorem~\ref{thm:hopf_main}).}
Simulation and numerical bifurcation analysis exhibit sustained oscillations
in the mixed mechanism~\cite{SuwanmajoKrishnan2015}, whose skeleton is the
transformation cycle $S_0\to S_2\to S_1\to S_0$ with the phosphatase on two of
its three reactions, contrary to the convergence half of the conjecture
of~\cite{ConradiShiu2015,ConradiShiu2018}.  Suwanmajo, Ramesh, and
Krishnan~\cite{SuwanmajoRameshKrishnan2020} found oscillations in the cyclic
network when a kinase or a phosphatase is shared, but not with distinct kinases and distinct
phosphatases.  Conradi, Mincheva, and Shiu~\cite{ConradiMinchevaShiu2019}
proved Hopf bifurcations for the mixed mechanism, and Conradi and
Mincheva~\cite{ConradiMincheva2024} for the cyclic network with a common kinase
and a common phosphatase after removing the dissociation steps; in both, supercriticality is observed in numerical
continuation rather than proved.
Theorem~\ref{thm:hopf_main} concerns the sharing pattern of the mixed
mechanism, realized with HMM modules; it proves that the Hopf bifurcation is supercritical, and
classifies the transformation cycles with at most three reactions.
Oscillatory regimes in multisite modules coupled to enzyme-activation
pathways~\cite{KrishnanSuwanmajo2018} involve mechanistic steps beyond the
cycles treated here.

\subsection{Organization}

The rest of the paper is organized as follows.
Section~\ref{sec:hypergraph} defines substrate hypergraphs, currents, and
cycles.  Section~\ref{sec:mech} defines HMM realizations and their
conservation laws, including the laws carried by directed cycles of the dual.
Section~\ref{sec:examples} gives examples of cycles and a two-reaction
non-cycle whose one-ray current cone nevertheless permits nonexistence of positive steady states and
multistationarity.  Section~\ref{sec:cycleproof} proves Theorem~\ref{thm:marked_main}.
Section~\ref{sec:convergence} proves
Theorem~\ref{thm:distinct_convergence} for cycles, and
Section~\ref{sec:instability} proves Theorem~\ref{thm:hopf_main} with
explicit rational rate constants and an interval-arithmetic enclosure of the first Lyapunov coefficient; simulation shows a stable periodic orbit (Figure~\ref{fig:hopf_cycle}).  Section~\ref{sec:extensions}
extends Theorems~\ref{thm:marked_main} and~\ref{thm:distinct_convergence} to weighted cycles, which allow higher stoichiometric coefficients, and to chains of
intermediates.  Section~\ref{sec:discussion} discusses open questions.

\section{Substrate hypergraphs, currents, and cycles}\label{sec:hypergraph}

\subsection{Substrate hypergraphs}
We recall the framework of substrate skeletons and hypergraph currents
introduced in~\cite{JoshiNguyen2026Hypergraph}, adapting the notation to the
present cycle setting.  The definitions are included to make the paper
self-contained; the dual graph and the cycle condition introduced below are
the structures specific to this work.

A \emph{directed hypergraph} $\mathcal H=(V,\mathcal E)$ consists of a finite
set of nodes $V$ and a finite set of directed hyperedges $\mathcal E$.  Each
hyperedge $\rho\in\mathcal E$ is an ordered pair
\[
  \rho=(t_\rho,h_\rho)
\]
of finite multisets of nodes, called its \emph{tail} and \emph{head}.  We write
$t_\rho(v)$ and $h_\rho(v)$ for the multiplicities of a node $v$ in these
multisets.  The arity of $\rho$ is $(|t_\rho|,|h_\rho|)$, where $|t_\rho|$
is the number of elements of $t_\rho$ counted with multiplicity.  We identify a
multiset with its formal sum; for example, $t_\rho=A+2B$ means that $A$ occurs
once and $B$ twice in the tail, and $|t_\rho|=3$.

Directed-hypergraph path and connectivity conventions are not uniform in the
literature; see, for example, Gallo et al.~\cite{GalloEtAl1993}.  We use the
following incidence-based convention.

\begin{definition}[Strong connectivity of a directed hypergraph]\label{def:strong-hypergraph}
A \emph{directed path} from $v_0$ to $v_k$ is a sequence
$v_0,\rho_1,v_1,\ldots,\rho_k,v_k$ with $v_{i-1}\in t_{\rho_i}$ and
$v_i\in h_{\rho_i}$ for all $i$.  The directed hypergraph is \emph{strongly
connected} if every ordered pair of nodes is joined by a directed path.
\end{definition}

Traversing a hyperedge requires entering through one tail node and leaving
through one head node; the other tail nodes play no role.  Equivalently, this
is reachability in the directed incidence graph obtained by inserting one
reaction-node for each hyperedge.  An ordinary directed graph is the special
case in which every edge has arity $(1,1)$; for ordinary directed graphs this
is the usual notion of strong connectivity.

A \emph{substrate hypergraph} encodes a set $\mathcal R$ of composite reactions
on a set $\mathcal S$ of substrate species.  The reactions are the hyperedges:
a reaction $\rho\in\mathcal R$ sends the nonempty reactant multiset $t_\rho$ to
the nonempty product multiset $h_\rho$,
\begin{equation}\label{eq:composite_reaction}
  \rho:\quad t_\rho\longrightarrow h_\rho.
\end{equation}
Thus $\mathcal H=(\mathcal S,\mathcal R)$ records which substrates each reaction
consumes and produces, without specifying its catalytic mechanism.  The
nonemptiness of tails and heads is part of the definition of a substrate
hypergraph; it excludes flow reactions $0\to A$, $A\to0$, and sequestration
reactions such as $A+B\to0$.  Throughout, $r\coloneqq |\mathcal R|$.

By arity, a reaction is a \emph{transformation} if $|t_\rho|=|h_\rho|=1$, a
\emph{transfer} if $|t_\rho|=|h_\rho|\ge2$, a \emph{binding} if
$|t_\rho|>|h_\rho|$, and a \emph{splitting} if $|t_\rho|<|h_\rho|$; a network
may combine all four.

The incidence matrix $M$ is the $|\mathcal S|\times r$ matrix whose rows are
indexed by substrates and whose columns are indexed by reactions:
\begin{equation}\label{eq:M_def}
  M_{X,\rho}=h_\rho(X)-t_\rho(X),
  \qquad X\in\mathcal S,\quad \rho\in\mathcal R.
\end{equation}
Thus the column indexed by $\rho$ is the net substrate reaction vector
$h_\rho-t_\rho$.

For a substrate $X$, define its total tail and head multiplicities by
\begin{equation}\label{eq:degrees}
  d_{\rm tail}(X)\coloneqq \sum_{\rho\in\mathcal R}t_\rho(X),
  \qquad
  d_{\rm head}(X)\coloneqq \sum_{\rho\in\mathcal R}h_\rho(X).
\end{equation}
These are multiset counts.  For example, in $2A\to B$, the substrate $A$
contributes two to $d_{\rm tail}(A)$.

The \emph{dual hypergraph} $\mathcal H^*$ has one node for each reaction
$\rho\in\mathcal R$ and one dual hyperedge for each substrate $X\in\mathcal S$.
The tail of the dual hyperedge $X$ contains every reaction that produces $X$,
with multiplicity $h_\rho(X)$, and its head contains every reaction that
consumes $X$, with multiplicity $t_\rho(X)$.  Thus a substrate runs in the dual
from its producers to its consumers, and the dual incidence matrix is
\begin{equation}\label{eq:dual_incidence_matrix}
  M^*=-M^{\mathsf T}.
\end{equation}
A dual hyperedge has an empty tail or head precisely when the corresponding
original node occurs in no head or no tail, respectively.  The dual of the dual
is $\mathcal H$.

The dual may have parallel hyperedges, for instance when two
substrates are both produced by one reaction and consumed by another, and
self-loops, when a substrate is produced and consumed by the same reaction.
When the dual is an ordinary directed graph, a \emph{directed cycle} is a
closed directed path with no repeated
reaction-node other than its initial and terminal node; a self-loop is a
directed cycle of length one.

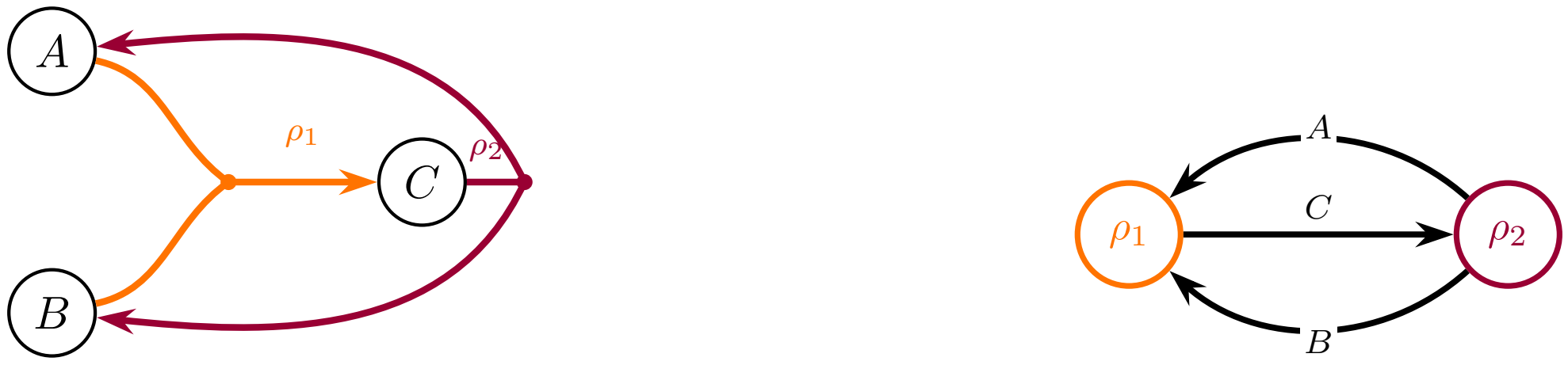
\begin{figure}[h!]
\centering
\begin{subfigure}[b]{0.53\textwidth}
\centering
\begin{tikzpicture}[
  nd/.style={circle,draw,thick,minimum size=0.72cm,inner sep=1pt},
  lab/.style={font=\scriptsize\itshape},
  scale=0.95
]
  \node[nd] (A) at (-2.7, 1.15) {$A$};
  \node[nd] (B) at (-2.7,-1.15) {$B$};
  \node[nd] (C) at ( 0.55,0) {$C$};
  \coordinate (m1) at (-1.15,0);
  \coordinate (m2) at ( 1.45,0);
  \draw[orange!90!red,line width=1.6pt] (A) to[out=-12,in=145] (m1);
  \draw[orange!90!red,line width=1.6pt] (B) to[out=12,in=215] (m1);
  \fill[orange!90!red] (m1) circle (2.0pt);
  \draw[-{Stealth[length=9pt,width=6pt]},orange!90!red,line width=1.6pt]
        (m1) -- (C) node[midway,above=4pt,lab]{$\rho_1$};
  \draw[purple!80!black,line width=1.6pt] (C) -- (m2);
  \fill[purple!80!black] (m2) circle (2.0pt);
  \draw[-{Stealth[length=9pt,width=6pt]},purple!80!black,line width=1.6pt]
        (m2) to[out=115,in=6] (A);
  \draw[-{Stealth[length=9pt,width=6pt]},purple!80!black,line width=1.6pt]
        (m2) to[out=245,in=-6] (B);
  \node[lab,purple!80!black] at (1.12,0.28) {$\rho_2$};
\end{tikzpicture}
\caption{The substrate hypergraph $\rho_1:A+B\to C$, $\rho_2:C\to A+B$.}
\end{subfigure}
\hfill
\begin{subfigure}[b]{0.41\textwidth}
\centering
\begin{tikzpicture}[
  rn/.style={circle,draw,line width=1.35pt,minimum size=0.86cm,inner sep=1pt,font=\small\itshape},
  el/.style={font=\scriptsize\itshape,fill=white,inner sep=1pt},
  scale=1.02
]
  \node[rn,draw=orange!90!red,text=orange!90!red] (r1) at (-1.55,0) {$\rho_1$};
  \node[rn,draw=purple!80!black,text=purple!80!black] (r2) at ( 1.55,0) {$\rho_2$};
  \draw[-{Stealth[length=9pt,width=6pt]},line width=1.45pt]
        (r1) -- node[el,above=2pt] {$C$} (r2);
  \draw[-{Stealth[length=9pt,width=6pt]},line width=1.45pt]
        (r2) to[bend left=42] (r1);
  \draw[-{Stealth[length=9pt,width=6pt]},line width=1.45pt]
        (r2) to[bend right=42] (r1);
  \node[el] at (0,0.88) {$A$};
  \node[el] at (0,-0.88) {$B$};
\end{tikzpicture}
\caption{The dual graph $\mathcal H^*$.}
\end{subfigure}
\caption{A substrate hypergraph and its dual.  Colors distinguish
hyperedges; branches joined at a dot are consumed or produced together.}
\label{fig:dual_convention}
\end{figure}

Figure~\ref{fig:dual_convention} shows the construction for a cycle containing
a binding and a splitting reaction.

\begin{lemma}[Strong connectivity and duality]\label{lem:strong-dual}
Let $\mathcal H$ be a directed hypergraph such that every hyperedge has
nonempty tail and nonempty head, and every node occurs in at least one tail and
at least one head.  Then $\mathcal H$ is strongly connected if and only if
$\mathcal H^*$ is strongly connected.
\end{lemma}
\begin{proof}
Suppose first that $\mathcal H$ is strongly connected.  Let $\rho$ and $\eta$
be hyperedges, which are nodes of $\mathcal H^*$.  Choose
$v\in h_\rho$ and $w\in t_\eta$.  Since $\mathcal H$ is strongly connected,
there is a directed path from $v$ to $w$.  Reading the nodes of that path as
dual hyperedges gives a directed path from reaction node $\rho$ to reaction
node $\eta$ in $\mathcal H^*$.

Conversely, suppose $\mathcal H^*$ is strongly connected.  By the hypotheses
the tail and head of the dual hyperedge corresponding to $v$ record,
respectively, the occurrences of $v$ in heads and tails of hyperedges of
$\mathcal H$.  Thus $\mathcal H^*$ satisfies the same nonemptiness
conditions.  Applying the first implication to $\mathcal H^*$
shows that $(\mathcal H^*)^*=\mathcal H$ is strongly connected.
\end{proof}

\subsection{Currents}
For each substrate, total production must balance total consumption at a
steady state.  This balance can be stated before any enzyme mechanism is
specified.

\begin{definition}[Current]\label{def:current}
A \emph{current} on $\mathcal H$ is a vector $J\in\mathbb R_{\ge0}^{r}$,
indexed by reactions, that satisfies \emph{flow conservation} at every
substrate $X$:
\begin{equation}\label{eq:flow_cons}
  \sum_{\rho\in\mathcal R}h_\rho(X)J_\rho
  =
  \sum_{\rho\in\mathcal R}t_\rho(X)J_\rho .
\end{equation}
Equivalently, $MJ=0$.  The currents form the convex polyhedral cone
\begin{equation}\label{eq:current_cone}
  \mathcal K(\mathcal H)
  \coloneqq \ker M\cap\mathbb R_{\ge0}^{r}.
\end{equation}
\end{definition}

Clarke's stoichiometric network analysis~\cite{Clarke1980} uses a
similar notion of current on the full reaction network; here currents are
defined on the substrate skeleton.  For the two-reaction futile skeleton in Figure~\ref{fig:intro_skeleton}, flow
conservation at either substrate gives $J_{\rho_1}=J_{\rho_2}$, so its current
cone is the single ray $\mathbb R_{\ge0}\one$.

When $\mathcal H^*$ is ordinary, $M^*$ is its node--edge incidence matrix,
with $-1$ at the tail and $+1$ at the head of each dual edge, and
by~\eqref{eq:dual_incidence_matrix}
\begin{equation}\label{eq:primal_dual_kernels}
  MJ=0 \iff (M^*)^{\mathsf T}J=0,
  \qquad
  M^{\mathsf T}p=0 \iff M^*p=0.
\end{equation}
Thus a current $J$ is a function on the reaction-nodes of the dual whose
values agree at the two ends of every dual edge, and a substrate conservation
law $p$ is a circulation on the dual edges: at every reaction-node, the total
of $p$ over incoming edges equals the total over outgoing edges.

\subsection{Cycles}\label{sec:cycles}
We use the following generalization of a directed cycle.

\begin{definition}[Cycle]\label{def:marked_cycle}
A substrate hypergraph $\mathcal H$ is a \emph{cycle} if it is strongly connected
and
\begin{equation}\label{eq:cycle_condition}
  d_{\rm tail}(X)=d_{\rm head}(X)=1
  \qquad(X\in\mathcal S).
\end{equation}
\end{definition}

Thus every substrate in a cycle is produced once and consumed once, with
multiplicity counted.  In particular, no substrate occurs with coefficient
greater than one in any tail or head, or in two different tails or two
different heads.  A reaction may consume or produce several different
substrates.  The futile skeleton of Figure~\ref{fig:intro_skeleton} and the
hypergraph of Figure~\ref{fig:dual_convention} are cycles; $2A\rightleftarrows B$
is not (Remark~\ref{rem:ray-not-cycle}).

\begin{proposition}[Cycle characterization]
\label{prop:cycle_char}
A substrate hypergraph $\mathcal H$ is a cycle if and only if its dual
$\mathcal H^*$ is a strongly connected ordinary directed graph.  In that case
\begin{equation}\label{eq:cycle_cone}
  \mathcal K(\mathcal H)=\mathbb R_{\ge0}\one,
\end{equation}
so every current is a nonnegative multiple of the all-ones current.
\end{proposition}
\begin{proof}
The dual hyperedge corresponding to $X$ has tail size $d_{\rm head}(X)$ and
head size $d_{\rm tail}(X)$.  Hence $\mathcal H^*$ is ordinary if and only if
$d_{\rm tail}(X)=d_{\rm head}(X)=1$ for every substrate $X$.  Strong
connectivity of $\mathcal H$ is then equivalent to strong connectivity of
$\mathcal H^*$ by Lemma~\ref{lem:strong-dual}.

By~\eqref{eq:primal_dual_kernels}, a current takes equal values at the two
ends of every dual edge.  Since the dual graph is connected, all entries of a
current are equal, so $J=\lambda\one$ with $\lambda\ge0$.  Conversely,
$\one$ satisfies flow conservation because every substrate has one
producer and one consumer.
\end{proof}

\begin{remark}\label{rem:ray-not-cycle}
The condition $\mathcal K(\mathcal H)=\mathbb R_{\ge0}\one$ does not
characterize cycles, even together with strong connectivity: $2A\to B$,
$B\to2A$ is strongly connected and has current cone
$\mathbb R_{\ge0}\one$ but $d_{\rm tail}(A)=d_{\rm head}(A)=2$.  A
one-dimensional current cone also does not imply the conclusion of
Theorem~\ref{thm:marked_main}; see Example~\ref{ex:ray_counterexample}.
\end{remark}

\section{Henri--Michaelis--Menten realizations}\label{sec:mech}

\subsection{Enzyme assignment and intermediates}\label{sec:enzyme_assign}
The detailed-model formalism in this section is the one developed
in~\cite{JoshiNguyen2026Hypergraph}, specialized to one-intermediate HMM
modules.  We recall it here while allowing arbitrary sharing of an enzyme
between composite reactions.

Let $E_1,\ldots,E_m$ be enzyme species.  An \emph{enzyme assignment} is a map
\begin{equation}\label{eq:sigma}
  \sigma:\mathcal R\longrightarrow\{0,1,\ldots,m\}
\end{equation}
whose image contains $\{1,\ldots,m\}$, so that every enzyme is used.  If
$\sigma(\rho)=a\ne0$, reaction $\rho$ is catalyzed by $E_a$; the value
$\sigma(\rho)=0$ means that the reaction is enzyme-free.  Different reactions
may share an enzyme.  For $a=1,\ldots,m$, write
\begin{equation}\label{eq:Lj}
  \mathcal R_a\coloneqq \{\rho\in\mathcal R:\sigma(\rho)=a\}\ne\varnothing
\end{equation}
for the reactions assigned to $E_a$.

Each reaction $\rho$ has its own intermediate $C_\rho$.  When $\rho$ is
enzyme-catalyzed, this intermediate is the bound form of the enzyme together
with the reactant multiset $t_\rho$; when $\rho$ is enzyme-free, it may be an
activated form of one or more substrates.  We reserve the word \emph{complex}
for its standard meaning in reaction network theory, a formal sum of species
such as $E+X$.  The substrates, the enzymes, and the intermediates are pairwise distinct
species.  Uppercase letters denote species and the corresponding lowercase
letters denote concentrations: a substrate $X$ has concentration $x_X$ (or $x$
in concrete examples), an enzyme $E_a$ has concentration $e_a$, and $C_\rho$
has concentration $c_\rho$.  For an enzyme-free reaction we write $E_0$ for
the absent enzyme, omit it from reaction schemes, and set $e_0\coloneqq 1$ in rate
formulas.

\subsection{Mass-action kinetics}
For a reaction $\rho\in\mathcal R$, the HMM module is
\begin{equation}\label{eq:gen_MM}
  \bar t_\rho\coloneqq E_{\sigma(\rho)}+t_\rho
  \;\xrightleftarrows{k_\rho^+}{k_\rho^-}\;
  C_\rho
  \;\xrightarrow{\,k_\rho\,}\;
  E_{\sigma(\rho)}+h_\rho\eqqcolon \bar h_\rho,
\end{equation}
with $\bar t_\rho=t_\rho$ and $\bar h_\rho=h_\rho$ when $\sigma(\rho)=0$.  The
rate constants $k_\rho^+,k_\rho^-,k_\rho$ are positive.

\begin{definition}[HMM realization]\label{def:realization}
Given a substrate hypergraph $\mathcal H=(\mathcal S,\mathcal R)$ and an enzyme
assignment $\sigma$, the \emph{HMM realization} $N_\sigma$ is the mass-action
network on the substrates, the enzymes $E_1,\ldots,E_m$, and the intermediates
$C_\rho$, $\rho\in\mathcal R$, with one module \eqref{eq:gen_MM} for every
$\rho\in\mathcal R$.  We write $\mathcal Z$ for its set of species, so that
concentration vectors lie in $\mathbb R_{\ge0}^{\mathcal Z}$.
We write $u\in\mathbb R_{\ge0}^{\mathcal Z}$ for such a vector; its
coordinates are the concentrations $x_X$, $e_a$, and $c_\rho$.
\end{definition}

Write $x=(x_X)_{X\in\mathcal S}$ for the vector of substrate
concentrations, the substrate part of $u$, and use the monomial notation
\begin{equation}\label{eq:monomial_notation}
  x^{t_\rho}\coloneqq \prod_{X\in\mathcal S}x_X^{t_\rho(X)}.
\end{equation}
The mass-action equation for the intermediate is
\begin{equation}\label{eq:ode_cl}
  \dot c_\rho
  =k_\rho^+e_{\sigma(\rho)}x^{t_\rho}
   -(k_\rho^-+k_\rho)c_\rho,
  \qquad \rho\in\mathcal R.
\end{equation}
We call
\begin{equation}\label{eq:kstar}
  J_\rho\coloneqq k_\rho c_\rho
  \qquad\text{and}\qquad
  k_\rho^*\coloneqq \frac{k_\rho k_\rho^+}{k_\rho^-+k_\rho}>0
\end{equation}
the \emph{catalytic flux} and \emph{composite catalytic rate} of reaction
$\rho$, respectively.  For a substrate $X$,
\begin{align}
  \dot x_X
  &=\sum_{\rho\in\mathcal R}h_\rho(X)k_\rho c_\rho
   +\sum_{\rho\in\mathcal R}t_\rho(X)k_\rho^-c_\rho \notag\\
  &\qquad
   -\sum_{\rho\in\mathcal R}t_\rho(X)k_\rho^+e_{\sigma(\rho)}x^{t_\rho}.
  \label{eq:ode_sS}
\end{align}
At a steady
state, \eqref{eq:ode_cl} gives
$k_\rho^+e_{\sigma(\rho)}x^{t_\rho}=(k_\rho^-+k_\rho)c_\rho$ for every
$\rho$.  Substituting this into~\eqref{eq:ode_sS} reduces the substrate
equations to
\begin{equation}\label{eq:ss_current_balance}
  MJ=0,
\end{equation}
where $J=(J_\rho)_{\rho\in\mathcal R}$.  Thus the catalytic flux vector at every
steady state is a current on $\mathcal H$.  At a steady state,
\eqref{eq:ode_cl} also gives
\begin{equation}\label{eq:throughput}
  J_\rho=k_\rho c_\rho
  =k_\rho^*e_{\sigma(\rho)}x^{t_\rho}.
\end{equation}
For the futile cycle, these formulas are
$J_{\rho_1}=k_1c_1=k_1^*ex$ and $J_{\rho_2}=k_2c_2=k_2^*fy$.

\subsection{Conservation laws}\label{sec:cons}
Every HMM realization has two families of linear conservation laws: enzyme
totals, and lifts of substrate conservation laws of the composite network.

Identify each species with a standard basis vector of $\mathbb R^{\mathcal Z}$
and each complex with the corresponding sum.  The stoichiometric subspace of
$N_\sigma$ is
\begin{equation}\label{eq:stoich_subspace}
  \mathcal W
  =\operatorname{span}\{\,C_\rho-\bar t_\rho,\ \bar h_\rho-C_\rho:
     \rho\in\mathcal R\,\}.
\end{equation}
A vector in $\mathcal W^\perp$ is a linear conservation law: its inner product
with the concentration vector is constant along every trajectory.  The value
of $L=\sum_Z\alpha_ZZ$ at a concentration vector is obtained by replacing each
species by its concentration.

\noindent\textbf{Enzyme conservation.}
For each enzyme $E_a$, the enzyme total
\begin{equation}\label{eq:gen_ej_cons}
  T_a\coloneqq E_a+\sum_{\rho\in\mathcal R_a}C_\rho,
  \qquad a=1,\ldots,m,
\end{equation}
is constant along trajectories.

\noindent\textbf{Composite and detailed substrate laws.}
A vector $p=(p_X)_{X\in\mathcal S}$ lies in $\ker M^{\mathsf T}$ exactly when its
weighted substrate total is unchanged by every composite reaction:
\begin{equation}\label{eq:composite_conservation_balance}
  \sum_X p_Xh_\rho(X)=\sum_X p_Xt_\rho(X)
  \qquad(\rho\in\mathcal R).
\end{equation}
For a vector $w\in\mathbb R^{\mathcal S}$ of substrate weights, define the
\emph{reactant weight} of reaction $\rho$ by
\begin{equation}\label{eq:tail_weight}
  \tau_\rho(w)\coloneqq \sum_{X\in\mathcal S}t_\rho(X)w_X .
\end{equation}
Thus $\tau_\rho(w)$ is the total $w$-weight of the reactant multiset, with
multiplicity counted: a tail $A+B$ has weight $w_A+w_B$, while $2A+B$ has
weight $2w_A+w_B$.  For $p\in\ker M^{\mathsf T}$, the number $\tau_\rho(p)$ is
the $p$-weight held in $C_\rho$ while the reactants are bound, and we define
the \emph{lifted law}
\begin{equation}\label{eq:lifted_law}
  L_p
  \coloneqq \sum_{X\in\mathcal S}p_XX
    +\sum_{\rho\in\mathcal R}\tau_\rho(p)C_\rho.
\end{equation}
If $p\ge0$, then every coefficient in $L_p$ is nonnegative.  For the futile
cycle, the substrate law $X+Y$ lifts to $X+Y+C_1+C_2$.

\begin{lemma}[All conservation laws in the HMM realization]\label{lem:lifted}
For every $p\in\ker M^{\mathsf T}$, the lifted law $L_p$ lies in
$\mathcal W^\perp$.  For every enzyme $E_a$, the enzyme total
$T_a$ of~\eqref{eq:gen_ej_cons} lies in $\mathcal W^\perp$.
Every linear conservation law of $N_\sigma$ has a unique representation
\begin{equation}\label{eq:all_conservation_laws}
  L_p+\sum_{a=1}^m q_aT_a,
  \qquad p\in\ker M^{\mathsf T},\quad q_a\in\mathbb R.
\end{equation}
\end{lemma}
\begin{proof}
Consider one HMM module for reaction $\rho$.  Binding removes substrate weight
$\tau_\rho(p)$ and creates $C_\rho$, whose coefficient in $L_p$ is
$\tau_\rho(p)$.  Thus binding, and hence unbinding, preserves $L_p$.  Catalysis
consumes $C_\rho$ and creates the product multiset $h_\rho$, whose weight is
$\sum_X p_X h_\rho(X)=\tau_\rho(p)$
by~\eqref{eq:composite_conservation_balance}.  So catalysis also preserves
$L_p$.  Binding replaces one free enzyme by one intermediate, unbinding reverses
that exchange, and catalysis replaces the intermediate by the free enzyme; hence
each $T_a$ is conserved.

Conversely, let a linear conservation law have substrate coefficients $p_X$,
enzyme coefficients $q_a$, and intermediate coefficients $\psi_\rho$, and set
$q_0\coloneqq 0$.  Orthogonality to the binding vector of reaction $\rho$ gives
$\psi_\rho=\tau_\rho(p)+q_{\sigma(\rho)}$.  Orthogonality to the catalytic vector
then gives~\eqref{eq:composite_conservation_balance}, so
$p\in\ker M^{\mathsf T}$.  The coefficients of the law therefore agree with
those of~\eqref{eq:all_conservation_laws}.  The substrate and free-enzyme
coordinates determine $p$ and $q$, so the representation is unique.
\end{proof}

\begin{corollary}[The conservation law carried by a dual cycle]
\label{cor:dcycle_law}
Let $\mathcal H$ be a cycle, so that its dual $\mathcal H^*$ is an
ordinary directed graph, and let $\gamma$ be a directed cycle in
the dual graph $\mathcal H^*$, with
substrate-edges $\mathcal S(\gamma)$ and reaction-nodes $\mathcal R(\gamma)$,
and let $p^\gamma\in\mathbb R^{\mathcal S}$ be the indicator vector of
$\mathcal S(\gamma)$.  Then $p^\gamma\in\ker M^{\mathsf T}$, and its lift is
the nonnegative conservation law
\begin{equation}\label{eq:dual_cycle_law}
  L_\gamma
  =\sum_{X\in\mathcal S(\gamma)}X
   +\sum_{\rho\in\mathcal R(\gamma)}C_\rho.
\end{equation}
\end{corollary}
\begin{proof}
Each reaction-node on $\gamma$ has exactly one incoming and one outgoing edge
of $\gamma$, and a reaction-node off $\gamma$ has none.  So $p^\gamma$ is a
circulation, and $M^{\mathsf T}p^\gamma=0$ by~\eqref{eq:primal_dual_kernels}.
The incoming dual edges at $\rho$ are the substrates of $t_\rho$, so
$\tau_\rho(p^\gamma)=1$ for $\rho\in\mathcal R(\gamma)$ and $0$ otherwise,
and~\eqref{eq:lifted_law} gives~\eqref{eq:dual_cycle_law}.
\end{proof}

For the futile cycle, the dual $2$-cycle gives $L_\gamma=X+Y+C_1+C_2$.

\subsection{Stoichiometric compatibility classes}\label{sec:scc}
For any $u^0\in\mathbb R_{\ge0}^{\mathcal Z}$, its \emph{stoichiometric
compatibility class} is
\begin{equation}\label{eq:scc}
  \mathcal C(u^0)
  \coloneqq (u^0+\mathcal W)\cap
    \mathbb R_{\ge0}^{\mathcal Z}.
\end{equation}
It is closed and convex.  We call it \emph{positive} if it contains a strictly
positive concentration vector.  Two nonnegative concentration vectors lie in
the same class exactly when their difference lies in $\mathcal W$.  Hence if
$L$ is a linear conservation law and $u^{(1)},u^{(2)}$ lie in the same class,
\begin{equation}\label{eq:compatibility_identity}
  L(u^{(1)})=L(u^{(2)}).
\end{equation}
In particular, every enzyme total and every nonzero nonnegative conservation
law takes a strictly positive value on a positive class, namely its value at a
strictly positive point of the class.

\begin{proposition}[Positive classes]\label{prop:positive_classes}
Let $\mathcal H$ be a cycle and let $u^0\in\mathbb R_{\ge0}^{\mathcal Z}$.
The compatibility class
$\mathcal C(u^0)$ is positive if and only if $T_a(u^0)>0$ for every enzyme
$E_a$ and $L_\gamma(u^0)>0$ for every directed cycle $\gamma$ of
$\mathcal H^*$.
\end{proposition}
\begin{proof}
If $\mathcal C(u^0)$ contains a strictly positive point $u$, then
by~\eqref{eq:compatibility_identity} each $T_a$ and each $L_\gamma$ takes at
$u^0$ its value at $u$, which is positive.

Conversely, suppose that $\mathcal C(u^0)$ is not positive.  Then the affine
set $u^0+\mathcal W$ does not meet the open orthant
$\mathbb R_{>0}^{\mathcal Z}$.  Separating these two convex sets gives a
nonzero vector $L$ with $L\cdot(u^0+w)\le L\cdot y$ for all $w\in\mathcal W$
and all $y>0$.  Letting $w$ range over $\mathcal W$ forces $L\in\mathcal
W^\perp$, letting one coordinate of $y$ grow forces $L\ge0$, and letting $y$
tend to $0$ gives $L(u^0)\le0$.  By Lemma~\ref{lem:lifted},
$L=L_p+\sum_aq_aT_a$ with $p\in\ker M^{\mathsf T}$.  The coefficients of $L$
on the substrates are the entries of $p$, and on $E_a$ the number $q_a$, so
$p\ge0$ and $q\ge0$.  By~\eqref{eq:primal_dual_kernels}, $p$ is a nonnegative
circulation on $\mathcal H^*$, and by the flow decomposition
theorem~\cite{AhujaMagnantiOrlin1993} it is a nonnegative combination
$\sum_\gamma\alpha_\gamma p^\gamma$ of indicator vectors of directed cycles.
Since the lift $p\mapsto L_p$ is linear,
$L=\sum_\gamma\alpha_\gamma L_\gamma+\sum_aq_aT_a$ with all coefficients
nonnegative and not all zero.  Then $L(u^0)>0$, a contradiction.
\end{proof}

\section{Examples}\label{sec:examples}

\begin{figure}[h!]
\centering
\begin{subfigure}[b]{0.43\textwidth}
\centering
\begin{tikzpicture}[baseline={(current bounding box.center)},scale=0.82]
  \node (S0) at (150:2.0cm) {$S_0$};
  \node (S1) at ( 30:2.0cm) {$S_1$};
  \node (S2) at (270:2.0cm) {$S_2$};
  \draw[-{Stealth[length=10pt,width=7pt]},line width=1.6pt,blue!80!black]
    (S0) .. controls +(0:1.2cm) and +(120:1.2cm) .. (S1)
    node[midway,above right=-2pt,font=\scriptsize\itshape,black] {$\rho_1$};
  \draw[-{Stealth[length=10pt,width=7pt]},line width=1.6pt,orange!85!red]
    (S1) .. controls +(240:1.2cm) and +(0:1.2cm) .. (S2)
    node[midway,right=2pt,font=\scriptsize\itshape,black] {$\rho_2$};
  \draw[-{Stealth[length=10pt,width=7pt]},line width=1.6pt,teal!80!black]
    (S2) .. controls +(120:1.2cm) and +(240:1.2cm) .. (S0)
    node[midway,left=2pt,font=\scriptsize\itshape,black] {$\rho_3$};
  \foreach \P in {S0,S1,S2}{\draw[thick] (\P) circle (0.32cm);}
\end{tikzpicture}
\caption{Transformation cycle, $N=3$ (Example~\ref{ex:trans_cycle}).}
\end{subfigure}
\hfill
\begin{subfigure}[b]{0.43\textwidth}
\centering
\begin{tikzpicture}[
  nd/.style={circle,draw,thick,minimum size=0.70cm,inner sep=1pt,font=\scriptsize},
  AH/.style={-{Stealth[length=8pt,width=5.5pt]},line width=1.45pt},
  LN/.style={line width=1.45pt},
  lbl/.style={font=\scriptsize\itshape,fill=white,inner sep=1pt},
  scale=0.73
]
  \node[nd] (X1p) at ( 90:2.5cm) {$X_{1,p}$};
  \node[nd] (X2p) at (-30:2.5cm) {$X_{2,p}$};
  \node[nd] (X3p) at (210:2.5cm) {$X_{3,p}$};
  \node[nd] (X1)  at ( 90:1.15cm) {$X_{1}$};
  \node[nd] (X2)  at (-30:1.15cm) {$X_{2}$};
  \node[nd] (X3)  at (210:1.15cm) {$X_{3}$};
  \coordinate (J1) at ( 30:1.85cm);
  \coordinate (J2) at (-90:1.85cm);
  \coordinate (J3) at (150:1.85cm);
  \draw[LN,orange!85!red] (X1p) -- (J1); \draw[LN,orange!85!red] (X2) -- (J1);
  \draw[AH,orange!85!red] (J1) -- (X1); \draw[AH,orange!85!red] (J1) -- (X2p);
  \fill[orange!85!red] (J1) circle (3pt);
  \node[lbl,text=orange!85!red] at (30:2.45cm) {$\rho_1$};
  \draw[LN,violet!80!black] (X2p) -- (J2); \draw[LN,violet!80!black] (X3) -- (J2);
  \draw[AH,violet!80!black] (J2) -- (X2); \draw[AH,violet!80!black] (J2) -- (X3p);
  \fill[violet!80!black] (J2) circle (3pt);
  \node[lbl,text=violet!80!black] at (-90:2.45cm) {$\rho_2$};
  \draw[LN,green!60!black] (X3p) -- (J3); \draw[LN,green!60!black] (X1) -- (J3);
  \draw[AH,green!60!black] (J3) -- (X3); \draw[AH,green!60!black] (J3) -- (X1p);
  \fill[green!60!black] (J3) circle (3pt);
  \node[lbl,text=green!60!black] at (150:2.45cm) {$\rho_3$};
\end{tikzpicture}
\caption{Transfer cycle, $n=3$ (Example~\ref{ex:transfer_cycle}).}
\end{subfigure}

\vspace{0.6em}
\begin{subfigure}[b]{0.72\textwidth}
\centering
\begin{tikzpicture}[
  nd/.style={circle,draw,minimum size=0.72cm,inner sep=1pt,font=\small,thick},
  c1/.style={orange!85!red},c2/.style={teal!80!black},c3/.style={violet!75!black},
  c4/.style={blue!65!black},c5/.style={green!55!black},
  leg1/.style={line width=1.6pt,orange!85!red},
  leg4/.style={line width=1.6pt,blue!65!black},
  arr1/.style={-{Stealth[length=9pt,width=6pt]},line width=1.6pt,orange!85!red},
  arr2/.style={-{Stealth[length=9pt,width=6pt]},line width=1.6pt,teal!80!black},
  arr3/.style={-{Stealth[length=9pt,width=6pt]},line width=1.6pt,violet!75!black},
  arr4/.style={-{Stealth[length=9pt,width=6pt]},line width=1.6pt,blue!65!black},
  arr5/.style={-{Stealth[length=9pt,width=6pt]},line width=1.6pt,green!55!black},
  junc/.style={circle,inner sep=0pt,minimum size=6pt},
  lbl/.style={font=\footnotesize\itshape},scale=0.82
]
  \node[nd] (Y4) at (-2.5, 1.0) {$Y_4$};
  \node[nd] (Y5) at (-2.5,-1.0) {$Y_5$};
  \node[nd] (Y2) at ( 0.0, 2.0) {$Y_2$};
  \node[nd] (Y1) at ( 0.0, 0.0) {$Y_1$};
  \node[nd] (Y3) at ( 0.0,-2.0) {$Y_3$};
  \node[nd] (X1) at ( 2.5, 1.0) {$X_1$};
  \node[nd] (X2) at ( 2.5,-1.0) {$X_2$};
  \coordinate (M1) at (0.833,1.0);
  \coordinate (S3) at (0.833,-1.0);
  \coordinate (M4) at (-1.3,0.0);
  \draw[leg1] (Y1) -- (M1); \draw[leg1] (Y2) -- (M1); \draw[arr1] (M1) -- (X1);
  \node[junc,fill=orange!85!red] at (M1) {};
  \node[lbl,c1] at ($(M1)!0.52!(X1)+(0,0.24)$) {$\rho_1$};
  \draw[arr2] (X1) -- (X2); \node[lbl,c2,right=3pt] at ($(X1)!0.5!(X2)$) {$\rho_2$};
  \draw[arr3] (X2) -- (S3); \draw[arr3] (S3) -- (Y1); \draw[arr3] (S3) -- (Y3);
  \node[junc,fill=violet!75!black] at (S3) {};
  \node[lbl,c3] at ($(X2)!0.5!(S3)+(0,0.25)$) {$\rho_3$};
  \draw[leg4] (Y3) -- (M4); \draw[leg4] (Y4) -- (M4);
  \draw[arr4] (M4) -- (Y2); \draw[arr4] (M4) -- (Y5);
  \node[junc,fill=blue!65!black] at (M4) {};
  \node[lbl,c4] at (-0.78,0.0) {$\rho_4$};
  \draw[arr5] (Y5) -- (Y4); \node[lbl,c5,left=3pt] at ($(Y5)!0.5!(Y4)$) {$\rho_5$};
\end{tikzpicture}
\caption{Mixed arity (Example~\ref{ex:mixed_cycle}).}
\end{subfigure}
\caption{Three cycles.  Colors distinguish reactions; branches joined at a dot
are consumed or produced together.}
\label{fig:three_families}
\end{figure}
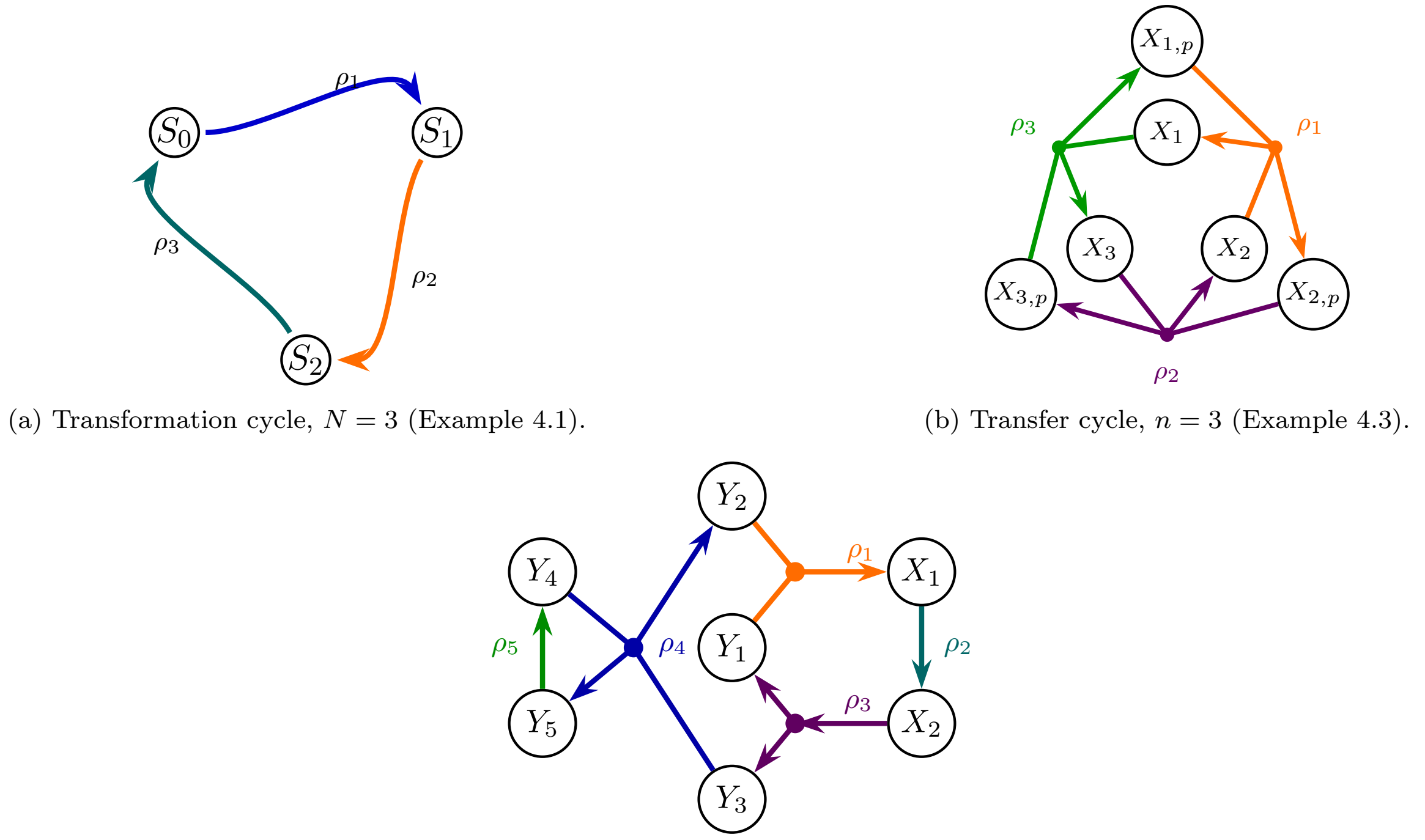

\begin{example}[Transformation cycle]\label{ex:trans_cycle}
On species $S_0,\dots,S_{N-1}$, with indices modulo $N$, let
\begin{equation}\label{eq:trans_cycle}
  \rho_\ell:S_{\ell-1}\longrightarrow S_\ell,
  \qquad \ell=1,\dots,N.
\end{equation}
Every substrate has one producer and one consumer, and the dual is a directed
$N$-cycle; Figure~\ref{fig:three_families}(a) shows the case $N=3$.
\exampleend
\end{example}

\begin{example}[Kinase--phosphatase pair with ATP regeneration]\label{ex:atp_cycle}
The reactions
\begin{equation}\label{eq:atp_cycle}
  \rho_1:X+\mathrm{ATP}\longrightarrow X_p+\mathrm{ADP},\qquad
  \rho_2:X_p\longrightarrow X,\qquad
  \rho_3:\mathrm{ADP}\longrightarrow \mathrm{ATP}
\end{equation}
form a cycle; see Figure~\ref{fig:intro_atp}.  Each of $X$, $X_p$, ATP, and ADP occurs once in a tail and once
in a head.  The dual graph has edges $X_p:\rho_1\to\rho_2$,
$X:\rho_2\to\rho_1$, $\mathrm{ADP}:\rho_1\to\rho_3$, and
$\mathrm{ATP}:\rho_3\to\rho_1$, so it is strongly connected.  Its two dual
cycles give the conserved totals $X+X_p+C_1+C_2$ and
$\mathrm{ATP}+\mathrm{ADP}+C_1+C_3$.  If ATP and ADP are held at fixed
concentrations, they are absorbed into the rate constants, and the skeleton
reduces to the futile cycle $X\rightleftarrows X_p$, which is again a cycle.
\exampleend
\end{example}

\begin{example}[Transfer cycle]\label{ex:transfer_cycle}
For $n\ge2$, on $2n$ species $X_1,\dots,X_n,X_{1,p},\dots,X_{n,p}$, consider
\begin{equation}\label{eq:transfer_cycle}
  \rho_\ell:X_{\ell,p}+X_{\ell+1}\longrightarrow X_\ell+X_{\ell+1,p},
  \qquad \ell=1,\dots,n,
\end{equation}
with indices modulo $n$.  Each substrate has one producer and one consumer.  In
the dual, $X_{\ell+1,p}$ is an edge $\rho_\ell\to\rho_{\ell+1}$ and $X_\ell$ is
an edge $\rho_\ell\to\rho_{\ell-1}$, so the dual is strongly connected; see
Figure~\ref{fig:three_families}(b).
\exampleend
\end{example}

\begin{example}[Mixed arity]\label{ex:mixed_cycle}
The reactions
\begin{equation}\label{eq:mixed_cycle}
  Y_1+Y_2\to X_1,\qquad
  X_1\to X_2,\qquad
  X_2\to Y_1+Y_3,\qquad
  Y_3+Y_4\to Y_2+Y_5,\qquad
  Y_5\to Y_4
\end{equation}
are a binding, a transformation, a splitting, a transfer, and a
transformation.  Every substrate has total tail and head multiplicity one, and
the dual graph is strongly connected; see Figure~\ref{fig:three_families}(c).
\exampleend
\end{example}

By Theorem~\ref{thm:marked_main}, every HMM realization of
Examples~\ref{ex:trans_cycle}--\ref{ex:mixed_cycle}, under any enzyme
assignment, has exactly one positive steady state in each positive
stoichiometric compatibility class and is persistent.

\begin{example}[Two futile cycles sharing both enzymes]\label{ex:feliu_wiuf}
The substrate hypergraph
\begin{equation}\label{eq:two_futile}
  X\to Y,\qquad Y\to X,\qquad Z\to W,\qquad W\to Z
\end{equation}
is not a cycle: its dual consists of two disjoint directed $2$-cycles, and its
current cone is two-dimensional.  If one enzyme catalyzes $X\to Y$ and
$Z\to W$ and another catalyzes $Y\to X$ and $W\to Z$, the HMM realization is
multistationary for some rate constants, as Feliu and Wiuf showed~\cite{FeliuWiuf2012}.  Sharing an
enzyme within a cycle, by contrast, never creates multistationarity, by
Theorem~\ref{thm:marked_main}.
\exampleend
\end{example}

\begin{example}[One-dimensional current does not imply unique existence]
\label{ex:ray_counterexample}
Consider the two-reaction substrate hypergraph
\begin{equation}\label{eq:ray_counterexample_skeleton}
  \rho_1:S_1+S_2\longrightarrow 2S_1,
  \qquad
  \rho_2:S_1+S_2\longrightarrow 2S_2.
\end{equation}
Flow conservation at either substrate gives $J_{\rho_1}=J_{\rho_2}$, so
$\mathcal K(\mathcal H)=\mathbb R_{\ge0}\one$.  The hypergraph is not a cycle:
each substrate has total tail multiplicity two and total head multiplicity two.

Realize the two reactions with distinct enzymes $E_1,E_2$, and write
$C_i\coloneqq C_{\rho_i}$, $c_i\coloneqq c_{\rho_i}$, $k_i\coloneqq k_{\rho_i}$, and
$k_i^*\coloneqq k_{\rho_i}^*$.  Take $k_1=k_2=1$ and write $T_i=e_i+c_i$ for the
two enzyme totals.  At a positive steady state the common current is a number
$J>0$, and~\eqref{eq:throughput} gives $c_1=c_2=J$ and
\begin{equation}\label{eq:ray_counterexample_flux0}
  J=k_1^*(T_1-J)s_1s_2,
  \qquad
  J=k_2^*(T_2-J)s_1s_2.
\end{equation}
Hence any positive steady state satisfies
\begin{equation}\label{eq:ray_counterexample_flux}
  k_1^*(T_1-J)=k_2^*(T_2-J)
\end{equation}
and
\begin{equation}\label{eq:ray_counterexample_product}
  s_1s_2=\frac{J}{k_1^*(T_1-J)}.
\end{equation}

\medskip
\noindent\emph{Failure of existence.}
Take $T_1=T_2=1$, $k_1^*=1$, and $k_2^*=2$.
Equation~\eqref{eq:ray_counterexample_flux} forces $J=1$, which gives
$e_1=e_2=0$.  Thus no positive steady state exists in any compatibility class
with these enzyme totals.

\medskip
\noindent\emph{Failure of uniqueness.}
Take $T_1=2$, $T_2=1$, $k_1^*=1$, and $k_2^*=10$.
Equation~\eqref{eq:ray_counterexample_flux} gives $J=8/9$, and hence
$c_1=c_2=8/9$, $e_1=10/9$, $e_2=1/9$, and $s_1s_2=4/5$.  On the compatibility
class where the conserved substrate total $S_1+S_2+2C_1+2C_2$ equals $50/9$,
we also have $s_1+s_2=2$.  Therefore the two substrate states are
\begin{equation}\label{eq:ray_two_states}
  (s_1,s_2)=\left(1+\frac1{\sqrt5},1-\frac1{\sqrt5}\right),
  \qquad
  \left(1-\frac1{\sqrt5},1+\frac1{\sqrt5}\right).
\end{equation}
The enzyme and intermediate coordinates are identical at the two states, and
their substrate difference is a multiple of
$S_1-S_2=h_{\rho_1}-t_{\rho_1}$, so the two positive steady states lie in the
same compatibility class.

These effective rates are realized by positive HMM rate constants: with
$k_i=k_i^-=1$, choosing $k_i^+=2k_i^*$ gives the specified values of $k_i^*$.
\exampleend
\end{example}

\section{Existence, uniqueness, and persistence}\label{sec:cycleproof}
Throughout this section $\mathcal H=(\mathcal S,\mathcal R)$ is a substrate
hypergraph cycle and $N_\sigma$ is an HMM realization of it.

\subsection{Uniqueness}
By~\eqref{eq:ss_current_balance} and Proposition~\ref{prop:cycle_char}, at every
steady state there is $\lambda\ge0$ with
\begin{equation}\label{eq:cycle_flux_identity}
  J_\rho=k_\rho c_\rho=k_\rho^*e_{\sigma(\rho)}x^{t_\rho}=\lambda
  \qquad(\rho\in\mathcal R).
\end{equation}

Let $u^{(1)}$ and $u^{(2)}$ be positive steady states, with common fluxes
$\lambda_1,\lambda_2>0$, and set
\begin{equation}\label{eq:theta_def}
  \theta\coloneqq \log\frac{\lambda_2}{\lambda_1}.
\end{equation}
For every species $Z$ with concentration $z$, define the logarithmic change
\begin{equation}\label{eq:mu_coordinate}
  \mu_Z\coloneqq \log\frac{z^{(2)}}{z^{(1)}},
\end{equation}
and set $\mu_{E_0}\coloneqq 0$.  Write $\mu=(\mu_X)_{X\in\mathcal S}$ for the substrate
part.  By~\eqref{eq:tail_weight},
\begin{equation}\label{eq:tail_log_change}
  \tau_\rho(\mu)
  =\sum_{X\in\mathcal S}t_\rho(X)\mu_X
  =\log\frac{(x^{(2)})^{t_\rho}}
                  {(x^{(1)})^{t_\rho}},
\end{equation}
so $\tau_\rho(\mu)>0$ means that the reactant monomial of $\rho$ is larger at
state $2$.

Taking logarithms of the ratios of~\eqref{eq:cycle_flux_identity} at the two
states gives, for every $\rho\in\mathcal R$,
\begin{equation}\label{eq:gen_relations}
  \mu_{C_\rho}=\theta,
  \qquad
  \mu_{E_{\sigma(\rho)}}+\tau_\rho(\mu)=\theta.
\end{equation}
Here and below $\mu_{E_0}=0$ for an enzyme-free reaction.

\begin{lemma}[Enzyme comparison]
\label{lem:tailpos}
Let $u^{(1)}$ and $u^{(2)}$ be positive steady states in the same
compatibility class with $\theta\ge0$.  Then $\mu_{E_a}\le0$ for every enzyme
$E_a$, with equality if and only if $\theta=0$.  Consequently
\begin{equation}\label{eq:tail_lower_bound}
  \tau_\rho(\mu)\ge\theta\ge0
  \qquad(\rho\in\mathcal R).
\end{equation}
\end{lemma}
\begin{proof}
The enzyme total~\eqref{eq:gen_ej_cons} has the same value at the two states,
and~\eqref{eq:gen_relations} gives
$c_\rho^{(2)}=\exp(\theta)c_\rho^{(1)}$.
Hence
\begin{equation}\label{eq:enzyme_comparison}
  e_a^{(2)}
  =e_a^{(1)}-(\exp(\theta)-1)
     \sum_{\rho\in\mathcal R_a}c_\rho^{(1)}.
\end{equation}
Since $\mathcal R_a\ne\varnothing$ and $c_\rho^{(1)}>0$, this gives
$e_a^{(2)}\le e_a^{(1)}$, with equality if and only if $\theta=0$.  The bound
\eqref{eq:tail_lower_bound} follows from~\eqref{eq:gen_relations}, since
$\mu_{E_{\sigma(\rho)}}\le0$, including $\mu_{E_0}=0$.
\end{proof}

\begin{proposition}[Uniqueness]\label{prop:unique}
Every HMM realization of a cycle has at most one positive steady state in each
stoichiometric compatibility class, for all positive rate constants and every
enzyme assignment.
\end{proposition}
\begin{proof}
Let $u^{(1)}$ and $u^{(2)}$ be positive steady states in the same compatibility
class.  Interchanging them changes the signs of $\theta$ and of every $\mu_Z$,
so we may assume $\theta\ge0$.  Lemma~\ref{lem:tailpos} gives
$\tau_\rho(\mu)\ge\theta\ge0$ for every $\rho$.

Suppose first that $\mu_X\ne0$ for some substrate $X$.  Give the dual edge $X$
the weight $\mu_X$.  The incoming dual edges at a reaction-node $\rho$ are the
substrates of $t_\rho$, each with multiplicity one, so the incoming sum at
$\rho$ is $\tau_\rho(\mu)\ge0$.  Lemma~\ref{lem:nonnegative_edge_cycle} of Appendix~\ref{app:graph} gives a
directed cycle $\gamma$ with $\mu_X\ge0$ for every $X\in\mathcal S(\gamma)$ and
$\mu_X>0$ for at least one.  Every intermediate in $L_\gamma$ has
$\mu_{C_\rho}=\theta\ge0$.  Thus every species in $L_\gamma$ is nondecreasing
and at least one increases, contradicting the equality of the conserved value
$L_\gamma$ at compatible states.

Hence $\mu_X=0$ for every substrate, so $\tau_\rho(\mu)=0$ for every $\rho$.
Then~\eqref{eq:tail_lower_bound} forces $\theta=0$, so
\eqref{eq:gen_relations} and Lemma~\ref{lem:tailpos} show that every
intermediate and every free enzyme is unchanged.  The two states are equal.
\end{proof}

The same comparison applies to two states with equal species-formation
rates, not only to two steady states.  Write $\dot u=f(u)$ for the
mass-action system of $N_\sigma$.

\begin{proposition}[Injectivity]\label{prop:injective}
Let $N_\sigma$ be an HMM realization of a cycle, under any enzyme assignment
and with any positive rate constants.  If $u^{(1)}$ and $u^{(2)}$ are strictly
positive, lie in the same stoichiometric compatibility class, and satisfy
$f(u^{(1)})=f(u^{(2)})$, then $u^{(1)}=u^{(2)}$.
\end{proposition}

Proposition~\ref{prop:unique} is the case $f(u^{(1)})=f(u^{(2)})=0$.

\begin{proof}
For a quantity $g$, write $\Delta g\coloneqq g(u^{(2)})-g(u^{(1)})$, and let
\[
  \phi_\rho(u)\coloneqq k_\rho^+e_{\sigma(\rho)}x^{t_\rho}-k_\rho^-c_\rho
\]
be the net binding rate of reaction $\rho$, with $e_0=1$.
By~\eqref{eq:ode_cl}, the component of $f$ at $C_\rho$ is
$\phi_\rho-J_\rho$, and by~\eqref{eq:ode_sS} the component at a substrate $X$
is $\sum_\rho h_\rho(X)J_\rho-\sum_\rho t_\rho(X)\phi_\rho$.  Equality of the
intermediate components gives $\Delta\phi_\rho=\Delta J_\rho$ for every
$\rho$, and the substrate components then give $M\,\Delta J=0$.
By~\eqref{eq:primal_dual_kernels} and the connectivity of the dual graph,
$\Delta J=\delta\one$ for some $\delta\in\mathbb R$.  Interchanging the two
states if necessary, we may assume $\delta\ge0$.  Then
\begin{equation}\label{eq:inject_intermediate}
  \Delta c_\rho=\frac{\delta}{k_\rho}\ge0
  \qquad(\rho\in\mathcal R),
\end{equation}
and the enzyme totals~\eqref{eq:gen_ej_cons}, which agree at the two states,
give $\Delta e_a=-\delta\sum_{\rho\in\mathcal R_a}k_\rho^{-1}\le0$.  The
identity $\Delta\phi_\rho=\delta$ reads
\begin{equation}\label{eq:inject_binding}
  k_\rho^+\Delta\bigl(e_{\sigma(\rho)}x^{t_\rho}\bigr)
  =\delta\Bigl(1+\frac{k_\rho^-}{k_\rho}\Bigr)\ge0.
\end{equation}
Since $e_{\sigma(\rho)}^{(2)}\le e_{\sigma(\rho)}^{(1)}$, this gives
$(x^{(2)})^{t_\rho}\ge(x^{(1)})^{t_\rho}$, that is, $\tau_\rho(\mu)\ge0$ with
$\mu_X=\log(x_X^{(2)}/x_X^{(1)})$ as in~\eqref{eq:mu_coordinate}.

If $\mu_X\ne0$ for some substrate $X$, the argument in the proof of
Proposition~\ref{prop:unique} gives a directed cycle $\gamma$ of the dual
graph on which every $\mu_X$ is nonnegative and at least one is positive.
By~\eqref{eq:inject_intermediate}, every intermediate in $L_\gamma$ is
nondecreasing, so $L_\gamma(u^{(2)})>L_\gamma(u^{(1)})$, contradicting
compatibility.  Hence $x^{(1)}=x^{(2)}$, and the left side
of~\eqref{eq:inject_binding} equals $k_\rho^+(x^{(1)})^{t_\rho}\Delta
e_{\sigma(\rho)}\le0$.  Therefore $\delta=0$,
and~\eqref{eq:inject_intermediate} and the enzyme totals give
$u^{(1)}=u^{(2)}$.
\end{proof}

\begin{remark}[Consequences of injectivity]\label{rem:injective}
By Proposition~\ref{prop:injective} and~\cite[Theorem~1.4]{MullerEtAl2016}
(see also~\cite{FeliuWiuf2012Injectivity}),
for all positive rate constants and every strictly positive state $u$, the
Jacobian $Df(u)$ has trivial kernel on $\mathcal W$.  Since $f$ takes values
in $\mathcal W$, $Df(u)$ maps $\mathcal W$ into itself, so its restriction to
$\mathcal W$ is invertible.  Three consequences follow.  First, every positive
steady state is nondegenerate: its Jacobian has no zero eigenvalue on
$\mathcal W$.  Second, as rate constants or conserved totals vary
continuously, no real eigenvalue on $\mathcal W$ crosses zero, so a positive
steady state can lose stability only when a pair of complex conjugate
eigenvalues crosses the imaginary axis; there are no saddle-node or
transcritical bifurcations of positive steady states.  The loss of stability
in Theorem~\ref{thm:hopf_main} is of this kind.  Third, by the implicit
function theorem, the positive steady state of a positive class depends
smoothly on the rate constants and the conserved totals.
\end{remark}

\subsection{Existence and boundary exclusion}
\begin{lemma}[A positive conservation law]\label{lem:marked_conservative}
There is a nonnegative conservation law $L^*$ whose coefficient is strictly
positive on every substrate and every intermediate.  Consequently
$L^*+\sum_{a=1}^mT_a$ is a conservation law with strictly positive coefficient
on every species.
\end{lemma}
\begin{proof}
Because $\mathcal H^*$ is strongly connected, every dual edge lies on a
directed cycle.  For each substrate $X$, choose a directed cycle $\gamma_X$
containing it and set $L^*\coloneqq \sum_{X\in\mathcal S}L_{\gamma_X}$.  Every
substrate has positive coefficient in $L^*$.  For a reaction $\rho$, choose
$X\in t_\rho$; the dual edge $X$ enters $\rho$, so $\gamma_X$ contains $\rho$
and $C_\rho$ has positive coefficient.  Every enzyme has coefficient one in
$\sum_aT_a$.
\end{proof}

The following standard fact is proved in Appendix~\ref{app:fixedpoint}.

\begin{lemma}[A compact convex invariant set contains an equilibrium]
\label{lem:compact_equilibrium}
Let $K\subset\mathbb R^n$ be a nonempty compact convex set that is positively
invariant under a continuously differentiable vector field $\dot u=f(u)$.
Then $K$ contains a point $u^*$ with $f(u^*)=0$.
\end{lemma}

\begin{proposition}[No boundary steady states in a positive class]
\label{prop:no_boundary}
Every steady state in a positive stoichiometric compatibility class is
strictly positive.
\end{proposition}
\begin{proof}
Let $u^*$ be such a steady state, with common flux $\lambda\ge0$ as
in~\eqref{eq:cycle_flux_identity}.

\smallskip
\noindent\emph{Case 1: $\lambda>0$.}
Then $c_\rho=\lambda/k_\rho>0$ and
$k_\rho^*e_{\sigma(\rho)}x^{t_\rho}=\lambda>0$ for every $\rho$.  Hence
every enzyme has positive free concentration, since every enzyme is used, and
every substrate in $t_\rho$ has positive concentration.  Every substrate
occurs in a tail, so every coordinate of $u^*$ is positive.

\smallskip
\noindent\emph{Case 2: $\lambda=0$.}
Then $c_\rho=0$ for every $\rho$.  The positive enzyme totals make every free
enzyme positive, so $x^{t_\rho}=0$ and every reactant multiset
contains a substrate of zero concentration.  The corresponding dual edges
enter every reaction-node, so Lemma~\ref{lem:backward_walk} of
Appendix~\ref{app:graph} gives a directed
cycle $\gamma$ consisting of zero-concentration substrates.  Every species in
$L_\gamma$ vanishes at $u^*$, so $L_\gamma(u^*)=0$.  But $L_\gamma$ is a
nonzero nonnegative conservation law, so it is positive on a positive class.
This contradiction rules out Case 2.
\end{proof}

\begin{proposition}[Existence]\label{prop:existence}
Every HMM realization of a cycle has at least one positive steady state in each
positive stoichiometric compatibility class, for all positive rate constants
and every enzyme assignment.
\end{proposition}
\begin{proof}
Fix a positive compatibility class $\mathcal C\coloneqq \mathcal C(u^0)$.  It is nonempty, closed,
and convex.  The fixed value of the strictly positive conservation law of
Lemma~\ref{lem:marked_conservative} bounds every nonnegative coordinate, so $\mathcal C$
is compact.  Mass-action dynamics preserve the nonnegative orthant and the
affine set $u^0+\mathcal W$, so $\mathcal C$ is positively invariant.
Lemma~\ref{lem:compact_equilibrium} gives a steady state $u^*\in \mathcal C$, and
Proposition~\ref{prop:no_boundary} shows that $u^*$ is strictly positive.
\end{proof}

\subsection{Persistence}\label{sec:persistence}
Persistence and boundary steady states in networks with intermediates and
catalysts are studied in~\cite{AngeliDeLeenheerSontag2007,ShiuSturmfels2010,
MarcondesDeFreitasFeliuWiuf2017}.  A reaction network is \emph{persistent} if, for every trajectory starting at a
strictly positive concentration vector, the $\omega$-limit set does not meet
the boundary of the nonnegative orthant~\cite{AngeliDeLeenheerSontag2007}.

\begin{lemma}[Boundary exclusion implies persistence]
\label{lem:compact-boundary-persistence}
Let $\mathcal C$ be a compact positive stoichiometric compatibility class of a
mass-action network.  If $\mathcal C$ contains no boundary steady state, then no
trajectory starting at a strictly positive point of $\mathcal C$ has an
$\omega$-limit point on the boundary of the nonnegative orthant.
\end{lemma}
\begin{proof}
Let $u(t)$ be such a trajectory.  Suppose its $\omega$-limit set contains a
boundary point $p$; then $p\in \mathcal C$, since $\mathcal C$ is closed and positively
invariant.  Let $\mathcal Z_0$ be the
nonempty set of species whose concentration vanishes at $p$.
Recall that a nonempty set of species is a \emph{siphon} if every elementary
reaction producing one of its species also consumes a species in the set.
By a result of Angeli, De~Leenheer, and Sontag~\cite[Proposition~5.4]{AngeliDeLeenheerSontag2007}, the zero set
$\mathcal Z_0$ of a boundary $\omega$-limit point is a siphon; see
also~\cite[Theorem~2.14]{AndersonShiu2010} for mass-action systems.  Consequently,
the face on which all species in $\mathcal Z_0$ vanish is positively invariant
\cite[Proposition~5.5]{AngeliDeLeenheerSontag2007}.

Let
\begin{equation}\label{eq:invariant_face}
  \mathcal F\coloneqq \{u\in \mathcal C:\ u_Z=0\ \text{for all } Z\in\mathcal Z_0\}.
\end{equation}
It is nonempty, since $p\in\mathcal F$, and it is compact, convex, and
positively invariant.  Lemma~\ref{lem:compact_equilibrium}
gives a steady state in $\mathcal F$.  This steady state lies on the boundary of the
class $\mathcal C$, contrary to the hypothesis.
\end{proof}

\begin{proposition}[Persistence]\label{prop:persistence}
Every HMM realization of a cycle is persistent, for all positive rate constants
and every enzyme assignment.
\end{proposition}
\begin{proof}
The compatibility class of every strictly positive initial condition is
compact by Lemma~\ref{lem:marked_conservative} and contains no boundary
steady state by Proposition~\ref{prop:no_boundary}.  The result follows
from Lemma~\ref{lem:compact-boundary-persistence}.
\end{proof}
For transformation cycles in which every reaction is catalyzed by an
enzyme, Propositions~\ref{prop:no_boundary}
and~\ref{prop:persistence} also follow from~\cite[Theorem~3.15]{PerezMillanDickenstein2018};
see Section~\ref{sec:related}.

\begin{lemma}[Instant positivity]\label{lem:instant_positivity}
Let $\mathcal C$ be a compact positive compatibility class of a mass-action network
that contains no boundary steady state.  Then every trajectory starting in $\mathcal C$
is strictly positive for all $t>0$.
\end{lemma}
\begin{proof}
For each species $Z$, the equation for its concentration has the form
$\dot z=P_Z(u)-zQ_Z(u)$, where $P_Z\ge0$ collects the production terms and
$Q_Z$ is a polynomial.  Hence once $z(t_0)>0$, also $z(t)>0$ for all
$t\ge t_0$, and the set $\mathcal Z(t)$ of species that vanish at time $t$ is
nonincreasing in $t$.  Fix $\varepsilon>0$.  Since $\mathcal Z(t)$ takes
finitely many values, it equals a fixed set $\mathcal Z_*$ on some open
interval $I\subset(0,\varepsilon)$.  For $Z\in\mathcal Z_*$, $z\equiv0$ on
$I$, so $P_Z(u(t))=0$ on $I$: every elementary reaction producing $Z$ has a
reactant in $\mathcal Z_*$.  Suppose $\mathcal Z_*\ne\varnothing$.  Then the
set $\mathcal Z_*$ is a siphon, so its zero face is positively invariant by
\cite[Proposition~5.5]{AngeliDeLeenheerSontag2007}.  Its intersection with
$\mathcal C$ is nonempty, compact, and convex, so
Lemma~\ref{lem:compact_equilibrium} gives a boundary steady state in
$\mathcal C$, a contradiction.  Thus $u(t)$ is strictly positive for $t\in I$.
To conclude, fix any $t>0$ and choose $\varepsilon<t$.  The interval
$I\subset(0,\varepsilon)$ contains some $t_0<t$ with
$\mathcal Z(t_0)=\varnothing$.  Since $\mathcal Z(\tau)$ is nonincreasing,
$\mathcal Z(t)=\varnothing$.  Hence the trajectory is strictly positive for
every $t>0$.
\end{proof}

\begin{proof}[Proof of Theorem~\ref{thm:marked_main}]
Proposition~\ref{prop:unique} gives at most one positive steady state in
each compatibility class, Proposition~\ref{prop:existence} gives at least one
in every positive class, Proposition~\ref{prop:no_boundary} excludes
boundary steady states in such a class, and Proposition~\ref{prop:persistence}
together with Lemma~\ref{lem:instant_positivity} gives part~2: a trajectory
starting in a positive class is strictly positive at every $t_0>0$, and
Proposition~\ref{prop:persistence} applies to the trajectory started at
$t_0$; Proposition~\ref{prop:injective} gives part~3.
\end{proof}

\begin{remark}\label{rem:duality}
The uniqueness proof and the boundary-exclusion argument, on which
existence and persistence rest, use the same combinatorial step.  Each of them selects a set of dual edges that enters every
reaction-node, or a weighting whose incoming sums are nonnegative, extracts a
directed dual cycle, and applies its conservation law $L_\gamma$.  The selected
edges are the substrates with nonnegative log-change in the uniqueness proof
and the zero substrates in the boundary argument.  Persistence and uniqueness
do not exclude periodic or other non-equilibrium recurrent behavior; the
instability example in Section~\ref{sec:instability} illustrates this distinction.
\end{remark}

\section{Global convergence with distinct enzymes}\label{sec:convergence}

Throughout this section $\mathcal H$ is a cycle, and the enzyme assignment has
\emph{distinct enzymes}: each enzyme $E_a$, $a\ge1$, catalyzes exactly one
reaction, while enzyme-free reactions are allowed.  We write $E_\rho$ for
$E_{\sigma(\rho)}$ and $e_\rho$ for its concentration when $\rho$ is
enzyme-catalyzed.  When $\rho$ is enzyme-free, every occurrence of $E_\rho$
below is omitted and we set $e_\rho\coloneqq 1$ in the rate formula.

The convergence argument is an application of the reaction-coordinate theory
of Angeli and Sontag~\cite{AngeliSontag2008} and Angeli, De~Leenheer, and
Sontag~\cite{AngeliDeLeenheerSontag2010}.  The former introduced the
translation-invariant monotone-system argument and applied it to the futile
cycle; the latter gave graph-theoretic conditions for monotonicity and a
general convergence theorem for persistent chemical reaction networks.  Our
contribution in this section is to show that the substrate-hypergraph cycle
axioms verify their kernel, monotonicity, and irreducibility hypotheses
uniformly for every distinct-enzyme HMM realization.  Lemma~\ref{lem:flux_kernel}
identifies the translation direction, while Lemma~\ref{lem:coop} derives the
required sign pattern and irreducibility directly from the dual graph.  For
transformation cycles in which each reaction has its own enzyme and a chain of
intermediates, Eithun and Shiu~\cite{EithunShiu2017} verified these hypotheses
in the same way, after removing intermediates
by~\cite{MarcondesDeFreitasWiufFeliu2017Convergence}; their result is a
special case of Corollary~\ref{cor:special_cases}.  For transformation
cycles, the stoichiometric matrix also has the factorization required by
Donnell and Banaji~\cite{DonnellBanaji2013} (Section~\ref{sec:related}).

\subsection{The futile cycle}\label{sec:futile_rc}

\begin{figure}[h!]
\centering
\tikzset{rx/.style={draw,rounded corners=2pt,minimum width=0.9cm,minimum height=0.6cm,font=\small},
  parr/.style={-{Stealth[length=6pt,width=4.5pt]},line width=0.9pt},
  narr/.style={-{Stealth[length=6pt,width=4.5pt]},line width=0.9pt,dashed,red!75!black},
  elab/.style={font=\scriptsize,fill=white,inner sep=1pt}}
\begin{subfigure}[b]{0.36\textwidth}
\centering
\begin{tikzpicture}
  \node[rx] (B1) at (0,1.6) {$B_1$};
  \node[rx] (K1) at (2.4,1.6) {$K_1$};
  \node[rx] (B2) at (2.4,0) {$B_2$};
  \node[rx] (K2) at (0,0) {$K_2$};
  \draw[parr] (B1) to[bend left=18] node[elab,above]{$C_1$} (K1);
  \draw[parr] (K1) to[bend left=18] node[elab,below]{$E,C_1$} (B1);
  \draw[parr] (B2) to[bend left=18] node[elab,below]{$C_2$} (K2);
  \draw[parr] (K2) to[bend left=18] node[elab,above]{$F,C_2$} (B2);
  \draw[parr] (K1) to[bend left=30] node[elab,right]{$Y$} (B2);
  \draw[parr] (K2) to[bend left=30] node[elab,left]{$X$} (B1);
\end{tikzpicture}
\caption{Futile cycle, distinct enzymes.}
\end{subfigure}\hfill
\begin{subfigure}[b]{0.6\textwidth}
\centering
\begin{tikzpicture}
  \node[rx] (B1) at (0,1.6) {$B_1$};
  \node[rx] (K1) at (2.2,1.6) {$K_1$};
  \node[rx] (B2) at (4.4,1.6) {$B_2$};
  \node[rx] (K2) at (6.6,1.6) {$K_2$};
  \node[rx] (B3) at (4.4,0) {$B_3$};
  \node[rx] (K3) at (2.2,0) {$K_3$};
  \draw[parr] (B1) to[bend left=18] (K1);
  \draw[parr] (K1) to[bend left=18] (B1);
  \draw[parr] (B2) to[bend left=18] (K2);
  \draw[parr] (K2) to[bend left=18] (B2);
  \draw[parr] (B3) to[bend left=18] (K3);
  \draw[parr] (K3) to[bend left=18] (B3);
  \draw[parr] (K1) to node[elab,above]{$X_2,E$} (B2);
  \draw[parr] (K2) to[bend left=14] node[elab,left]{$X_3$} (B3);
  \draw[parr] (K3) to[bend left=14] node[elab,left]{$X_1$} (B1);
  \draw[parr] (K2) to[out=110,in=70,looseness=0.55] node[elab,above]{$E$} (B1);
  \draw[narr,{Stealth[length=6pt,width=4.5pt]}-{Stealth[length=6pt,width=4.5pt]}] (B1) to[out=-25,in=-155,looseness=0.8] node[elab,below]{$E$} (B2);
\end{tikzpicture}
\caption{Three-cycle with $E$ shared by $\rho_1,\rho_2$.}
\end{subfigure}
\caption{Sign graphs of the Jacobian $DF$ in reaction coordinates.  An
arrow $k\to j$, for $k\ne j$, means $\partial F_j/\partial\xi_k\ne0$; where shown,
the label gives the species through which elementary reaction $k$ affects the rate of $j$.  Solid
arrows are positive and dashed arrows negative.  (a)~With distinct enzymes
every arrow is positive, and the graph is strongly connected
(Lemma~\ref{lem:coop}).  (b)~In the network~\eqref{eq:unstable-network},
where $E$ catalyzes $\rho_1$ and $\rho_2$ and $F$ catalyzes $\rho_3$, the
binding steps $B_1$ and $B_2$ compete for $E$, which produces the negative
arrows in both directions between them, drawn as one double-headed arrow.}
\label{fig:influence}
\end{figure}
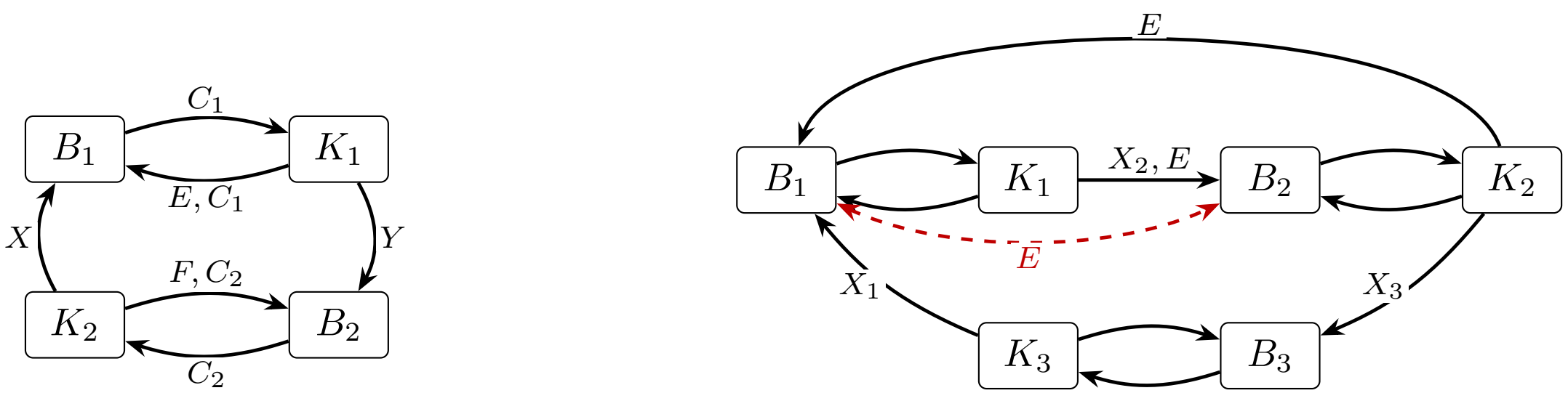

Angeli and Sontag treated the futile cycle in reaction coordinates
in~\cite[Section~3]{AngeliSontag2008}.  Splitting each module into its binding
and catalytic steps gives
\begin{equation}\label{eq:futile_BK}
\begin{aligned}
  &B_1:\ E+X\rightleftarrows C_1,\qquad &&K_1:\ C_1\to E+Y,\\
  &B_2:\ F+Y\rightleftarrows C_2,\qquad &&K_2:\ C_2\to F+X,
\end{aligned}
\end{equation}
and, with species ordered $(X,Y,E,F,C_1,C_2)$ and steps ordered
$(B_1,K_1,B_2,K_2)$, the stoichiometric
matrix is
\begin{equation}\label{eq:futile_Gamma}
  \Gamma=\begin{pmatrix}
  -1&0&0&1\\
  0&1&-1&0\\
  -1&1&0&0\\
  0&0&-1&1\\
  1&-1&0&0\\
  0&0&1&-1
  \end{pmatrix}.
\end{equation}
Here $\Gamma\one=0$ for $\one=(1,1,1,1)$, and the reaction-coordinate
Jacobian has the positive, strongly connected sign graph in
Figure~\ref{fig:influence}(a).  Lemma~\ref{lem:flux_kernel} generalizes the
kernel identity, and Lemma~\ref{lem:coop} generalizes the sign computation to
every distinct-enzyme realization of a substrate cycle.

\subsection{Reaction coordinates}
Split each module~\eqref{eq:gen_MM} into its binding reaction and its
catalytic reaction,
\begin{equation}\label{eq:BK}
  B_\rho:\ E_\rho+t_\rho\rightleftarrows C_\rho,
  \qquad
  K_\rho:\ C_\rho\to E_\rho+h_\rho ,
\end{equation}
with net rates
\begin{equation}\label{eq:BK_rates}
  R_{B_\rho}(u)=k_\rho^+e_\rho\,x^{t_\rho}-k_\rho^-c_\rho,
  \qquad
  R_{K_\rho}(u)=k_\rho c_\rho .
\end{equation}
These are the net binding rate $\phi_\rho$ of
Proposition~\ref{prop:injective} and the catalytic flux $J_\rho$
of~\eqref{eq:kstar}; we write them as $R_{B_\rho}$ and $R_{K_\rho}$ to follow
the reaction-coordinate notation of~\cite{AngeliDeLeenheerSontag2010}.
For an enzyme-free module, display~\eqref{eq:BK} means
$B_\rho:t_\rho\rightleftarrows C_\rho$ and
$K_\rho:C_\rho\to h_\rho$.
Let $\mathcal J\coloneqq \{B_\rho,K_\rho:\rho\in\mathcal R\}$, and let
$\Gamma\in\mathbb R^{\mathcal Z\times\mathcal J}$ be the matrix whose column
for $B_\rho$ is $C_\rho-E_\rho-t_\rho$ and whose column for $K_\rho$ is
$E_\rho+h_\rho-C_\rho$.  The HMM realization is
\begin{equation}\label{eq:ode_Gamma}
  \dot u=\Gamma R(u).
\end{equation}
For the futile cycle, $\Gamma$ is~\eqref{eq:futile_Gamma}.

\begin{lemma}[Flux kernel]\label{lem:flux_kernel}
$\ker\Gamma=\mathbb R\one$, where $\one\in\mathbb R^{\mathcal J}$ is the vector with every
entry equal to one.  At a positive steady state $u^*$, $R(u^*)=\lambda\one$ for
some $\lambda>0$.
\end{lemma}
\begin{proof}
If $\Gamma J=0$, the row of $C_\rho$ gives $J_{B_\rho}=J_{K_\rho}\eqqcolon J_\rho$,
and the rows of the substrates give $MJ=0$.  By~\eqref{eq:primal_dual_kernels}, $J$ takes equal
values at the two ends of every dual edge, and the dual graph is connected, so
$J\in\mathbb R\one$.  At $u^*$, $\Gamma R(u^*)=0$ gives $R(u^*)=\lambda\one$, and
$\lambda=k_\rho c^*_\rho>0$.
\end{proof}
Lemma~\ref{lem:flux_kernel} is the elementary-reaction counterpart of
Proposition~\ref{prop:cycle_char}.  A vector $J\in\ker\Gamma$ has
$J_{B_\rho}=J_{K_\rho}$ for every $\rho$, and the common values satisfy flow
conservation~\eqref{eq:flow_cons}.  Thus the nonnegative vectors in
$\ker\Gamma$ are the currents on $\mathcal H$, each value $J_\rho$ repeated on
the binding and catalytic steps of $\rho$, and $\one\in\mathbb R^{\mathcal J}$
corresponds to the current $\one\in\mathbb R^{\mathcal R}$.  The statement $R(u^*)=\lambda\one$ refines the fact from
Section~\ref{sec:mech} that the catalytic fluxes at a steady state form a
current: at a steady state the net binding rate of each module also equals its
catalytic flux.

The state $u(t)$ records the current concentrations but not how far each
elementary reaction has run.  Reaction coordinates record the latter.  Fix a
strictly positive initial condition $u^0$, and let $u(t)$ be the solution
of~\eqref{eq:ode_Gamma} with $u(0)=u^0$.  For each elementary reaction
$j\in\mathcal J$, let
\begin{equation}\label{eq:xi_def}
  \xi_j(t)\coloneqq \int_0^tR_j(u(s))\,ds
\end{equation}
be the net amount by which $j$ has proceeded up to time $t$; in chemistry this
is its extent of reaction.  Following~\cite{AngeliSontag2008}, we call
$\xi=(\xi_j)_{j\in\mathcal J}$ the \emph{reaction coordinates}.

The reaction coordinates determine the state.  Each elementary reaction $j$
moves $u$ along its reaction vector $\Gamma_{\cdot,j}$, by the amount
$\xi_j(t)$ up to time $t$, so integrating~\eqref{eq:ode_Gamma} gives
\begin{equation}\label{eq:u_from_xi}
  u(t)=u^0+\Gamma\xi(t).
\end{equation}
The reaction coordinates also
satisfy a closed system: differentiating~\eqref{eq:xi_def} and
using~\eqref{eq:u_from_xi} gives
\begin{equation}\label{eq:xi_system}
  \dot\xi=F(\xi)\coloneqq R(u^0+\Gamma\xi),\qquad \xi(0)=0,
\end{equation}
on the convex set
\[
X_{u^0}\coloneqq \{\xi\in\mathbb R^{\mathcal J}:\
u^0+\Gamma\xi\in\mathbb R_{\ge0}^{\mathcal Z}\},
\qquad
D_{u^0}\coloneqq \operatorname{int}X_{u^0}.
\]
The map $\xi\mapsto u^0+\Gamma\xi$ is not injective.  Adding $s\one$ to $\xi$
advances every elementary reaction by the same amount $s$; each module then
converts $s$ units of $t_\rho$ into $h_\rho$ and returns its enzyme and
intermediate to their former values, and each substrate is produced once and
consumed once, so the state does not change; that is, $\Gamma\one=0$.  Since $F$
depends on $\xi$ only through $u^0+\Gamma\xi$,
\begin{equation}\label{eq:translation_invariance}
  F(\xi+s\one)=F(\xi)\qquad(s\in\mathbb R).
\end{equation}
Monotonicity is applied to~\eqref{eq:xi_system}, not to~\eqref{eq:ode_Gamma}:
Lemma~\ref{lem:coop} shows that the Jacobian of $F$ has nonnegative
off-diagonal entries and is irreducible.

\begin{lemma}[Cooperativity and irreducibility]\label{lem:coop}
For every $\xi\in X_{u^0}$, every off-diagonal entry of $DF(\xi)$ is
nonnegative.  For every $\xi\in D_{u^0}$, the matrix $DF(\xi)$ is
irreducible.
\end{lemma}
\begin{proof}
Let $u=u^0+\Gamma\xi$, which is nonnegative.  The entry of $DF(\xi)$ in
row $j$ and column $k$ is
\begin{equation}\label{eq:coop_entry}
  \frac{\partial F_j}{\partial\xi_k}
  =\sum_{Z\in\mathcal Z}\frac{\partial R_j}{\partial u_Z}(u)\,\Gamma_{Z,k}.
\end{equation}
We first compute the row of a catalytic reaction $K_\rho$.
$R_{K_\rho}$ depends only on $c_\rho$, which is changed
only by $B_\rho$ and $K_\rho$.  The only nonzero off-diagonal entry is
$\partial F_{K_\rho}/\partial\xi_{B_\rho}=k_\rho>0$.

Next we compute the row of a binding reaction $B_\rho$, one column at a
time.  $R_{B_\rho}$ increases with $e_\rho$ and with $x_X$ for
$X\in t_\rho$, and decreases with $c_\rho$.
\begin{itemize}
\item Column $K_\rho$: $K_\rho$ produces $E_\rho$ and $h_\rho$ and consumes
$C_\rho$, so every term of~\eqref{eq:coop_entry} is nonnegative, and the
$C_\rho$ term equals $k_\rho^->0$.
\item Column $B_{\rho'}$, $\rho'\ne\rho$: $B_{\rho'}$ changes only $E_{\rho'}$,
$C_{\rho'}$, and the substrates of $t_{\rho'}$.  None of these is an argument of
$R_{B_\rho}$, because $E_{\rho'}\ne E_\rho$ by distinct enzymes and
$t_{\rho'}\cap t_\rho=\varnothing$ because each substrate lies in one tail.
The entry is zero.
\item Column $K_{\rho'}$, $\rho'\ne\rho$: $K_{\rho'}$ changes $E_{\rho'}$,
$C_{\rho'}$, and the substrates of $h_{\rho'}$, which it produces.  The only
arguments of $R_{B_\rho}$ among these are the substrates in
$h_{\rho'}\cap t_\rho$, each contributing a nonnegative term.  Thus all
off-diagonal entries are nonnegative.  When $u$ is strictly positive, the
entry is positive exactly when $\rho'\to\rho$ is an edge of the dual graph.
\end{itemize}
Now suppose $\xi\in D_{u^0}$, so that $u$ is strictly positive, and draw an
arrow $k\to j$ whenever $\partial F_j/\partial\xi_k>0$.  The arrows
include $B_\rho\to K_\rho$ and $K_\rho\to B_\rho$ for every $\rho$, and
$K_{\rho'}\to B_\rho$ for every dual edge $\rho'\to\rho$.  Since the dual
graph is strongly connected, so is this directed graph.  Figure~\ref{fig:influence}(a)
shows it for the futile cycle.
\end{proof}

Thus Lemma~\ref{lem:coop} directly verifies, for the present family, the
conclusion of the graph-theoretic monotonicity criterion
in~\cite[Corollary~1]{AngeliDeLeenheerSontag2010}.

The sign rule behind Lemma~\ref{lem:coop} is that a species shared by two
elementary reactions contributes a nonnegative term to~\eqref{eq:coop_entry}
when it lies on opposite sides of them, and a nonpositive term when it lies on
the same side.  With distinct enzymes, every species lies in exactly two
elementary reactions, on opposite sides: $C_\rho$ in $B_\rho$ and $K_\rho$;
$E_\rho$ in $B_\rho$ and $K_\rho$; a substrate $X$ in
$B_{\rho_{\rm cons}(X)}$ and $K_{\rho_{\rm prod}(X)}$.

Write $\xi\le\eta$ if $\xi_j\le\eta_j$ for all $j$, and $\xi\ll\eta$ if
$\xi_j<\eta_j$ for all $j$.   Solutions of~\eqref{eq:xi_system} exist for all $t\ge0$, because
$\dot\xi=R(u(t))$ is bounded on the compact compatibility class, and
solutions starting in $D_{u^0}$ remain there, because strictly positive
concentrations remain strictly positive.  A cooperative vector field on an
open convex set generates a monotone flow, and a strongly monotone one if its
Jacobian is also irreducible~\cite[Theorems~3.12 and~3.13]{HirschSmith2005}.
By Lemma~\ref{lem:coop}, the reaction-coordinate system is therefore strongly
monotone on $D_{u^0}$.
Thus, if $\xi(t)$ and $\eta(t)$ solve~\eqref{eq:xi_system} in $D_{u^0}$ with
$\xi(0)\le\eta(0)$ and $\xi(0)\ne\eta(0)$, then $\xi(t)\ll\eta(t)$ for all
$t>0$.  It is also monotone on $X_{u^0}$.  If $\xi\le\eta$ in $X_{u^0}$ and
$\zeta\in D_{u^0}$, then $(1-\epsilon)\xi+\epsilon\zeta\le(1-\epsilon)\eta+\epsilon\zeta$
for $0<\epsilon\le1$, and both points lie in $D_{u^0}$ because they correspond
to strictly positive concentrations.  Monotonicity on $D_{u^0}$ and continuous
dependence on initial conditions, as $\epsilon\to0$, give monotonicity on
$X_{u^0}$.

\subsection{Application of the convergence theorem}

\begin{proof}[Proof of Theorem~\ref{thm:distinct_convergence}]
By Lemma~\ref{lem:instant_positivity}, a trajectory starting in a
positive class is strictly positive at every positive time, so it suffices to
treat strictly positive initial conditions.  For such an initial condition,
by Lemma~\ref{lem:coop} and the paragraph following it, the reaction-coordinate system is monotone on
$X_{u^0}$ and strongly monotone on $D_{u^0}$.  Moreover,
$\ker\Gamma=\mathbb R\one$ by Lemma~\ref{lem:flux_kernel}, and
$\one$ lies in the interior of the positive orthant.  Thus the kernel
condition~(9) of~\cite{AngeliDeLeenheerSontag2010},
$\ker\Gamma\cap\operatorname{int}\mathbb R^{\mathcal J}_{\ge0}\ne\varnothing$, holds.  The original system
is persistent by Theorem~\ref{thm:marked_main}.

Angeli, De~Leenheer, and Sontag~\cite[Section~2]{AngeliDeLeenheerSontag2010}
fix the reactant side and the product side of each reversible reaction and
assume:
\begin{itemize}
\item[(H1)] no species is both a reactant and a product of the same reaction;
\item[(H2)] each rate is $C^1$; an irreversible rate is nondecreasing in each
  of its reactants and independent of all other species, and a reversible rate
  is a forward rate of this kind in its reactants minus a backward rate of this
  kind in its products;
\item[(H3)] an irreversible or forward rate vanishes when some reactant
  concentration is zero, and a backward rate vanishes when some product
  concentration is zero;
\item[(H4)] every stoichiometric compatibility class is compact;
\item[(H5)] the partial derivatives in (H2) are positive at every strictly
  positive state.
\end{itemize}
Under (H1)--(H5), \cite[Theorem~2, case~(9)]{AngeliDeLeenheerSontag2010}
states the following.  Suppose the network is persistent, and its
reaction-coordinate system is monotone on $X_{u^0}$ and strongly monotone on
the interior of $X_{u^0}$, with respect to the order of an orthant $\mathcal O$.
Suppose also that $\ker\Gamma\cap\operatorname{int}\mathcal O\ne\varnothing$.  Then
every solution starting at a strictly positive state converges to an
equilibrium, and this equilibrium is unique in its compatibility class.

Regard each $B_\rho$ as one reversible reaction, with forward rate
$k_\rho^+e_\rho x^{t_\rho}$ and backward rate $k_\rho^-c_\rho$, and each
$K_\rho$ as an irreversible reaction.  Then (H1) holds because the tail and
head of each elementary reaction have disjoint supports, (H2), (H3), and (H5)
hold for mass-action rates, and (H4) holds by
Lemma~\ref{lem:marked_conservative}.  The remaining hypotheses are the
monotonicity, kernel, and persistence statements verified above, so the
theorem applies.
\end{proof}

\begin{remark}[The underlying contraction argument]
Theorem~\ref{thm:distinct_convergence} rests on the
translation-invariant monotone-system theorem of Angeli and
Sontag~\cite[Theorem~1]{AngeliSontag2008}.  For the nonnegative orthant and
translation direction $\one$, their gauge gives the oscillation seminorm
\[
  \operatorname{osc}(\zeta)\coloneqq \max_j\zeta_j-\min_j\zeta_j,
  \qquad \Gamma\zeta=u-u^*.
\]
It is well defined because $\zeta$ is unique modulo $\one$, and their
contraction argument makes it strictly decrease away from the translation
direction~\cite[Lemmas~2.2--2.3]{AngeliSontag2008}.  We use the later general
convergence theorem rather than repeat that argument.

Theorems~13 and~14 of Al-Radhawi and Angeli~\cite{AlRadhawiAngeli2016}, on
piecewise-linear Lyapunov functions (see also~\cite{BlanchiniGiordano2014}),
give another direct route.  If the binding steps are first regarded as
irreversible, Lemma~\ref{lem:flux_kernel} gives a one-dimensional positive
flux kernel and every species is consumed by a unique elementary reaction.
Each binding step has the sole product $C_\rho$ and is its unique producing
reaction, so Theorem~14 permits adding its reverse; strong connectivity of
the dual gives the required common-ancestor condition.  Their
piecewise-linear-in-rates Lyapunov function, together with the positive
conservation law of Lemma~\ref{lem:marked_conservative} and the unique steady
state of Theorem~\ref{thm:marked_main}, therefore gives the same
conclusion.
\end{remark}

When two reactions $\rho,\rho'$ share an enzyme $E_a$, the argument fails;
Figure~\ref{fig:influence}(b) shows the resulting sign graph for the network
of Section~\ref{sec:instability}.
Both binding reactions consume $E_a$, so
$\partial F_{B_\rho}/\partial\xi_{B_{\rho'}}=-k_\rho^+x^{t_\rho}<0$,
while $\partial F_{B_\rho}/\partial\xi_{K_{\rho'}}>0$ and
$\partial F_{K_{\rho'}}/\partial\xi_{B_{\rho'}}>0$.  These three signs rule
out cooperativity with respect to every orthant of $\mathbb R^{\mathcal J}$.
Section~\ref{sec:instability} shows that this is not only a limitation of the
method.

\section{Instability and Hopf bifurcation under enzyme sharing}\label{sec:instability}

The conclusions of Theorem~\ref{thm:marked_main} are uniform over enzyme
assignments, but stability is not.  The following example gives an exact
counterexample to global convergence for the full class.  It complements the
oscillations computed in~\cite{SuwanmajoRameshKrishnan2020} and the Hopf
bifurcations proved for the mixed mechanism~\cite{ConradiMinchevaShiu2019} and
for the cyclic network after removing the dissociation
steps~\cite{ConradiMincheva2024}: the example below uses only HMM modules, has
rational data and an exact instability certificate, and the bifurcation is
proved to be supercritical.

\begin{proposition}[A shared-enzyme cycle with an unstable steady state]
\label{prop:unstable-shared}
Consider the substrate transformation cycle
$X_1\to X_2\to X_3\to X_1$ with HMM realization
\begin{equation}\label{eq:unstable-network}
\begin{aligned}
 E+X_1&\;\xrightleftarrows{k_1^+}{k_1^-}\;C_1
       \;\xrightarrow{\,k_1\,}\;E+X_2,\\
 E+X_2&\;\xrightleftarrows{k_2^+}{k_2^-}\;C_2
       \;\xrightarrow{\,k_2\,}\;E+X_3,\\
 F+X_3&\;\xrightleftarrows{k_3^+}{k_3^-}\;C_3
       \;\xrightarrow{\,k_3\,}\;F+X_1.
\end{aligned}
\end{equation}
For the rate constants
\begin{equation}\label{eq:unstable-rates}
 (k_1^+,k_1^-,k_1)=\left(2400,\frac12,\frac52\right),\quad
 (k_2^+,k_2^-,k_2)=\left(\frac{102}{5},2,100\right),\quad
 (k_3^+,k_3^-,k_3)=\left(\frac{11}{100},\frac12,5\right),
\end{equation}
the compatibility class with conserved totals
\begin{equation}\label{eq:unstable-totals}
 x_1+x_2+x_3+c_1+c_2+c_3=91,\qquad
 e+c_1+c_2=\frac{411}{20},\qquad f+c_3=60
\end{equation}
has the unique positive steady state
\begin{equation}\label{eq:unstable-state}
 (x_1,x_2,x_3,e,f,c_1,c_2,c_3)
 =\left(\frac12,50,10,\frac1{20},50,20,\frac12,10\right),
\end{equation}
and this steady state is unstable within its compatibility class.
Consequently, not every positive trajectory of an HMM realization of a
substrate hypergraph cycle converges to its positive steady state.
\end{proposition}

\begin{proof}
At the point~\eqref{eq:unstable-state}, the net binding flux and the catalytic
flux both equal $50$ in every module.  Indeed, for the three modules these
pairs of fluxes are
\[
 2400\frac1{20}\frac12-\frac12(20)=\frac52(20)=50,
\]
\[
 \frac{102}{5}\frac1{20}(50)-2\frac12=100\frac12=50,
 \qquad
 \frac{11}{100}(50)(10)-\frac12(10)=5(10)=50.
\]
Thus~\eqref{eq:unstable-state} is a steady state, and the three equalities
in~\eqref{eq:unstable-totals} hold.  Uniqueness follows from
Theorem~\ref{thm:marked_main}.

To determine stability, eliminate $x_3,e,f$ using~\eqref{eq:unstable-totals}
and use $(x_1,x_2,c_1,c_2,c_3)$ as coordinates on the compatibility class.
The reduced Jacobian at~\eqref{eq:unstable-state} is
\begin{equation}\label{eq:unstable-jacobian}
A=\begin{pmatrix}
-120&0&2401/2&1200&5\\
0&-51/50&2045/2&1022&0\\
120&0&-1203&-1200&0\\
0&51/50&-1020&-1122&0\\
-11/2&-11/2&-11/2&-11/2&-121/10
\end{pmatrix}.
\end{equation}
Its characteristic polynomial is
\begin{align}\label{eq:unstable-characteristic}
\det(\nu I-A)={}&\nu^5+\frac{61453}{25}\nu^4
+\frac{145890901}{500}\nu^3+\frac{464285352}{125}\nu^2\notag\\
&+\frac{259727307}{20}\nu+350030340.
\end{align}
The signs in the first column of its Routh table, from the $\nu^5$ row to the
constant row, are
\[
 +,\ +,\ +,\ +,\ -,\ +.
\]
There are therefore exactly two roots in the open right half-plane.  For
reference, they are numerically
$1.3144436556\pm8.9284847378\,\mathrm{i}$; the other three roots have negative
real part.  Hence the steady state is a hyperbolic saddle-focus on its
five-dimensional compatibility class and is unstable.  Its stable manifold
has codimension two, so positive points arbitrarily close to it fail to
converge to it.
\end{proof}

\begin{proposition}[Hopf bifurcation]\label{prop:hopf}
In~\eqref{eq:unstable-network}, keep $(k_1^+,k_1^-,k_1)$, $(k_2^+,k_2^-,k_2)$,
and $k_3$ as in~\eqref{eq:unstable-rates}, and set
\begin{equation}\label{eq:hopf-family}
 k_3^-=s,\qquad k_3^+=\frac{s+5}{50}\qquad(s>0).
\end{equation}
The value $s=1/2$ gives~\eqref{eq:unstable-rates}.  For every $s>0$, the
point~\eqref{eq:unstable-state} is the unique positive steady state of the
class~\eqref{eq:unstable-totals}.  There is exactly one $s^*>0$,
$s^*=10.4673375245\ldots$, at which this steady state changes stability: it is
unstable for $0<s<s^*$ and asymptotically stable for $s>s^*$.  At $s=s^*$ the
Jacobian on the compatibility class has a simple pair of eigenvalues
$\pm i\omega$, $\omega\approx10.2755$, its other three eigenvalues have
negative real part, and the pair crosses the imaginary axis with nonzero
speed.  The first Lyapunov coefficient at $s=s^*$ is negative.  Consequently the
bifurcation is supercritical: for every $s<s^*$ sufficiently close to $s^*$,
the class~\eqref{eq:unstable-totals} contains a unique limit cycle near the
steady state, of amplitude of order $(s^*-s)^{1/2}$, and this limit cycle is
orbitally asymptotically stable within the class.
\end{proposition}

\begin{proof}
For every $s$, the binding flux of the third module at~\eqref{eq:unstable-state}
is $\frac{s+5}{50}(50)(10)-10s=50$, so~\eqref{eq:unstable-state} remains a
steady state with the same totals, and it is the unique positive one by
Theorem~\ref{thm:marked_main}.  In the coordinates $(x_1,x_2,c_1,c_2,c_3)$ of
the proof of Proposition~\ref{prop:unstable-shared}, the characteristic
polynomial is $\nu^5+a_1\nu^4+a_2\nu^3+a_3\nu^2+a_4\nu+a_5$.
Appendix~\ref{app:dynamics_computation} records the exact coefficients and
Hurwitz determinants.  They show that, for $s>0$, all $a_i$ and
$\Delta_1,\Delta_2,\Delta_3$ are positive, while
$\Delta_4=(9/31250000)g(s)$ and $\Delta_5=a_5\Delta_4$, where $g$ is the
quartic in~\eqref{eq:hopf-q}.  Its coefficients change sign exactly once, so by
Descartes' rule of signs $g$ has exactly one positive root $s^*$, counted
with multiplicity.  Hence $s^*$ is a simple root, $g'(s^*)\ne0$, $g<0$ on
$(0,s^*)$, and $g>0$ on $(s^*,\infty)$.  The first column of the Routh table
is $1,\Delta_1,\Delta_2/\Delta_1,\Delta_3/\Delta_2,\Delta_4/\Delta_3,a_5$.  For
$s>s^*$ all six entries are positive, so every eigenvalue has negative real
part.  For $0<s<s^*$ the signs are $+,+,+,+,-,+$, so exactly two eigenvalues
have positive real part~\cite{Gantmacher1959}.  The conditions $a_5(s^*)>0$,
$\Delta_1(s^*),\Delta_2(s^*),\Delta_3(s^*)>0$, $\Delta_4(s^*)=0$, and
$\Delta_4'(s^*)\ne0$ are the hypotheses of Liu's criterion~\cite{Liu1994},
which gives the stated eigenvalue structure at $s^*$ and the transversal
crossing.

Let $A$ be the Jacobian on the compatibility class at $s=s^*$, let $B(u,v)$ be
the symmetric bilinear form given by the second derivatives of the reduced
vector field at the steady state, and write $\langle u,v\rangle=\sum_i\bar
u_iv_i$.  Choose $q,p\in\mathbb C^5$ with $Aq=i\omega q$,
$A^{\mathsf T}p=-i\omega p$, $\langle q,q\rangle=1$, and $\langle
p,q\rangle=1$.  The reduced vector field is quadratic, so its third
derivatives vanish, and the first Lyapunov coefficient~\cite{Kuznetsov2004} is
\begin{equation}\label{eq:lyapunov}
 \ell_1=\frac{1}{2\omega}\,\mathrm{Re}\Bigl(
 \bigl\langle p,B\bigl(\bar q,(2i\omega I-A)^{-1}B(q,q)\bigr)\bigr\rangle
 -2\bigl\langle p,B\bigl(q,A^{-1}B(q,\bar q)\bigr)\bigr\rangle\Bigr).
\end{equation}
Appendix~\ref{app:dynamics_computation} encloses $s^*$ in a rational interval
of width less than $10^{-60}$ and evaluates every quantity
in~\eqref{eq:lyapunov} in midpoint--radius interval
arithmetic~\cite{Johansson2017}.  The resulting enclosure of $\ell_1$ has
center $-2.68605139416826643\ldots\times10^{-3}$ and radius less than
$10^{-22}$, so $\ell_1<0$.  The eigenvalue pair $\pm i\omega$ crosses from
the right half-plane into the left half-plane as $s$ increases through $s^*$,
and $\ell_1<0$, so the Hopf bifurcation
theorem~\cite{MarsdenMcCracken1976,Kuznetsov2004} gives, for $s<s^*$ close to
$s^*$, a unique limit cycle near the steady state, of amplitude of order
$(s^*-s)^{1/2}$, that is orbitally asymptotically stable on the
two-dimensional center manifold.  The other three eigenvalues have negative
real part, so the center manifold attracts nearby trajectories in the class,
and the limit cycle is orbitally asymptotically stable within the class.
\end{proof}

\begin{figure}[h!]
\centering
\includegraphics[width=\textwidth]{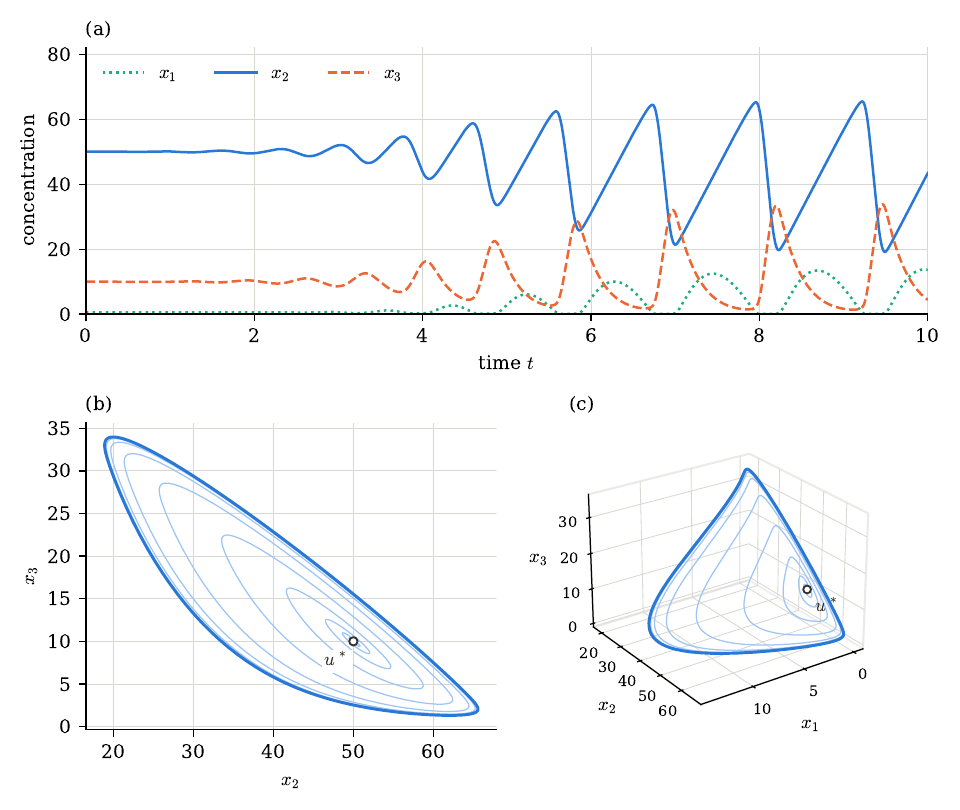}
\caption{Numerical integration of the shared-enzyme three-cycle~\eqref{eq:unstable-network}
with the rate constants~\eqref{eq:unstable-rates} ($s=1/2$), started from the steady
state~\eqref{eq:unstable-state} with $x_1$ increased and $x_2$ decreased by
$10^{-2}$, which leaves the totals~\eqref{eq:unstable-totals} unchanged.  (a)~Substrate concentrations for
$0\le t\le 10$; the oscillation grows and settles onto a periodic orbit of
period approximately $1.28$.  (b)~Projection onto the $(x_2,x_3)$ plane: the
transient (light) spirals out from the steady state $u^*$ (open circle) to the
periodic orbit (dark), shown for $30\le t\le 40$.  (c)~The same trajectory in
the substrate space $(x_1,x_2,x_3)$.}
\label{fig:hopf_cycle}
\end{figure}

Direct integration is consistent with Proposition~\ref{prop:hopf}: the
amplitude scales as $(s^*-s)^{1/2}$ near the crossing, and at $s=1/2$
trajectories near~\eqref{eq:unstable-state} approach a periodic orbit of period
approximately $1.28$ along which $x_2$ ranges from about $19$ to about $66$
(Figure~\ref{fig:hopf_cycle}).  Proposition~\ref{prop:hopf} concerns $s$ near
$s^*$; the periodic orbit at $s=1/2$ is computed, not proved.
Appendix~\ref{app:dynamics_computation} describes the reproducible
calculations.

The second enzyme in~\eqref{eq:unstable-network} is not needed.  Replace
the third module by the enzyme-free module
\begin{equation}\label{eq:free-module}
 X_3\;\xrightleftarrows{k_3^+}{k_3^-}\;C_3\;\xrightarrow{\,k_3\,}\;X_1,
\end{equation}
whose binding rate $k_3^+x_3-k_3^-c_3$ is first order in $x_3$.  The free enzyme
$F$ and its total are removed, and $x_3$ and $e$ are eliminated using the
first two totals in~\eqref{eq:unstable-totals}.

\begin{proposition}[The enzyme-free variant]\label{prop:hopf-free}
In the network obtained from~\eqref{eq:unstable-network} by replacing its third
module with~\eqref{eq:free-module}, keep $(k_1^+,k_1^-,k_1)$,
$(k_2^+,k_2^-,k_2)$, and $k_3=5$ as in~\eqref{eq:unstable-rates}, and set
\begin{equation}\label{eq:hopf-free-family}
 k_3^-=s,\qquad k_3^+=s+5\qquad(s>0).
\end{equation}
For every $s>0$, the point
$(x_1,x_2,x_3,e,c_1,c_2,c_3)=(\frac12,50,10,\frac1{20},20,\frac12,10)$ is the
unique positive steady state of the class with $x_1+x_2+x_3+c_1+c_2+c_3=91$
and $e+c_1+c_2=\frac{411}{20}$.  There is exactly one $\tilde s^*>0$,
$\tilde s^*=13.0160798148\ldots$, at which this steady state changes
stability: it is unstable for $0<s<\tilde s^*$, in particular at $s=1/2$, and
asymptotically stable for $s>\tilde s^*$.  At $s=\tilde s^*$ the Jacobian on
the compatibility class has a simple pair of eigenvalues $\pm i\tilde\omega$,
$\tilde\omega\approx10.7872$, its other three eigenvalues have negative real
part, the pair crosses the imaginary axis with nonzero speed, and the first
Lyapunov coefficient is negative.  Consequently the conclusions of
Proposition~\ref{prop:hopf} hold with $s^*$ replaced by $\tilde s^*$.
\end{proposition}

\begin{proof}
The binding flux of the third module at the stated point is
$(s+5)(10)-10s=50$, equal to the catalytic flux $5(10)$, and the first two
modules are unchanged, so the point is a steady state with the stated totals;
it is the unique positive one by Theorem~\ref{thm:marked_main}.  The proof of
Proposition~\ref{prop:hopf} applies verbatim to the data recorded in
Appendix~\ref{app:dynamics_computation}: all $a_i$ and
$\Delta_1,\Delta_2,\Delta_3$ are positive for $s>0$, $\Delta_4$ is a positive
multiple of a quartic $\tilde g$ whose coefficients change sign exactly once,
and the interval computation encloses $\ell_1$ in an interval with center
$-2.94698791053829051\ldots\times10^{-3}$ and radius less than $10^{-22}$.
\end{proof}

At $s=1/2$, direct integration shows trajectories near the steady state
approaching a periodic orbit of period approximately $1.33$ along which $x_2$
ranges from about $16$ to about $66$.  As in the two-enzyme case, the
periodic orbit at $s=1/2$ is computed, not proved.

\begin{proposition}[One enzyme on every reaction]\label{prop:convergent}
Let $\mathcal H$ be a transformation cycle with $n\ge1$ reactions, realized by
HMM modules in which one enzyme catalyzes all $n$ reactions.  Then for all
positive rate constants every trajectory in a positive compatibility class
converges to the positive steady state of that class.
\end{proposition}

\begin{proof}
The realization is the single ring of reactions
\[
 E+X_1\rightleftarrows C_1\to E+X_2\rightleftarrows C_2\to\cdots\to E+X_n
 \rightleftarrows C_n\to E+X_1,
\]
which is weakly reversible with one linkage class and $2n$ complexes.  The
differences between $E+X_1$ and the other $2n-1$ complexes, namely
$C_i-E-X_1$ for $1\le i\le n$ and $X_i-X_1$ for $2\le i\le n$, are linearly
independent, since each $C_i$ occurs in only one of them.  Hence the
stoichiometric subspace has dimension $2n-1$, and the deficiency is
$2n-1-(2n-1)=0$: the number of complexes, minus the number of linkage
classes, minus the dimension of the stoichiometric subspace.  By the
deficiency zero theorem~\cite{Feinberg1987} the system is complex balanced
for all positive rate constants, so the Horn--Jackson
function~\cite{HornJackson1972} is a Lyapunov function on each positive
class.  Its derivative along trajectories vanishes in the interior of a class
only at the complex-balanced steady state.  The class is compact, and by
persistence the $\omega$-limit set of a trajectory lies in its interior, so by
LaSalle's invariance principle it consists of that steady state.
\end{proof}

\begin{proposition}[The remaining assignment patterns]\label{prop:remaining}
Let $\mathcal H$ be a transformation cycle with at most three reactions, and
realize it by HMM modules under any enzyme assignment, with enzyme-free
reactions allowed.  Unless $\mathcal H$ has three reactions and some enzyme
catalyzes exactly two of them, for all positive rate constants every
trajectory in a positive compatibility class converges to the positive steady
state of that class.
\end{proposition}

\begin{proof}
Writing $E_1,E_2,E_3$ for distinct enzymes and $0$ for an enzyme-free
reaction, the assignments in which no enzyme catalyzes more than one reaction
are $(E_1)$ and $(0)$ on one reaction; $(E_1,E_2)$, $(E_1,0)$, and $(0,0)$ on
two reactions; and $(E_1,E_2,E_3)$, $(E_1,E_2,0)$, $(E_1,0,0)$, and $(0,0,0)$
on three reactions, in every order.  The remaining assignments are
$(E_1,E_1)$, $(E_1,E_1,E_1)$, and, up to rotation, the excluded patterns
$(E_1,E_1,E_2)$ and $(E_1,E_1,0)$.  The assignments in the first list are
covered by Theorem~\ref{thm:distinct_convergence}, and $(E_1,E_1)$ and
$(E_1,E_1,E_1)$ by Proposition~\ref{prop:convergent}.
\end{proof}

Three reactions is therefore the smallest size of a transformation cycle at
which convergence can fail.

\begin{proof}[Proof of Theorem~\ref{thm:hopf_main}]
Relabeling the substrates cyclically, we may assume that the enzyme
catalyzing two reactions catalyzes $X_1\to X_2$ and $X_2\to X_3$.  If the
third reaction has its own enzyme, parts~1 and~2 are
Propositions~\ref{prop:unstable-shared} and~\ref{prop:hopf}; if it is
enzyme-free, they are Proposition~\ref{prop:hopf-free}.  In both cases the
limit cycle of part~2 is a trajectory that does not converge to the steady
state.  This gives the ``if'' direction of part~3, since a three-reaction
transformation cycle in which an enzyme catalyzes exactly two reactions is,
up to rotation, one of these two networks.  The ``only if'' direction is
Proposition~\ref{prop:remaining}.
\end{proof}

\section{Extensions}\label{sec:extensions}

\subsection{Higher stoichiometric coefficients}
For a substrate $X$ that occurs in exactly one tail and exactly one head, let
$\rho_{\rm prod}(X)$ and $\rho_{\rm cons}(X)$ be the reactions producing and
consuming it, and write
\begin{equation}\label{eq:weights_ab}
  a_X\coloneqq h_{\rho_{\rm prod}(X)}(X),
  \qquad
  b_X\coloneqq t_{\rho_{\rm cons}(X)}(X)
\end{equation}
for its stoichiometric coefficients there.

\begin{definition}[Weighted cycle]\label{def:weighted_cycle}
A substrate hypergraph is a \emph{weighted cycle} if it is strongly connected,
every substrate occurs in exactly one reaction tail and exactly one reaction
head, with any positive coefficients, and $\ker M$ contains a strictly
positive vector $v$.
\end{definition}

A cycle is a weighted cycle with $a_X=b_X=1$ and $v=\one$.  For a
weighted cycle, let $\mathcal H^*$ now denote the directed graph with one node
per reaction and one edge $X$ from $\rho_{\rm prod}(X)$ to
$\rho_{\rm cons}(X)$ per substrate; it is the dual with the multiplicities
$a_X,b_X$ forgotten, and it is strongly connected by
Lemma~\ref{lem:strong-dual}.  Flow conservation at $X$ reads
\begin{equation}\label{eq:weighted_flow}
  a_X J_{\rho_{\rm prod}(X)}=b_X J_{\rho_{\rm cons}(X)},
\end{equation}
so a current is determined by its value at one reaction, and
$\mathcal K(\mathcal H)=\mathbb R_{\ge0}v$.

\begin{example}[A weighted cycle that is not a cycle]\label{ex:2A_B}
The substrate hypergraph $\rho_1:2A\to B$, $\rho_2:B\to2A$ of
Remark~\ref{rem:ray-not-cycle} is a weighted cycle, with $a_A=b_A=2$,
$a_B=b_B=1$, and $v=\one$.  Its dual is the $2$-cycle with edge $A$ from
$\rho_2$ to $\rho_1$ and edge $B$ from $\rho_1$ to $\rho_2$, and flow
conservation~\eqref{eq:weighted_flow} at $A$ and at $B$ reads
$2J_{\rho_2}=2J_{\rho_1}$ and $J_{\rho_1}=J_{\rho_2}$.  By contrast,
$\rho_1:A\to2B$, $\rho_2:B\to A$ has every substrate in exactly one tail and
one head, but flow conservation at $B$ gives $2J_{\rho_1}=J_{\rho_2}$ and at
$A$ gives $J_{\rho_2}=J_{\rho_1}$, so the only current is zero and the
hypergraph is not a weighted cycle.  The hypergraph
of Example~\ref{ex:ray_counterexample} is not a weighted cycle, since each
substrate occurs in two tails.
\exampleend
\end{example}

\begin{lemma}[Positive-current criterion]\label{lem:weighted_current_criterion}
The kernel of $M$ contains a strictly positive vector if and only if
\begin{equation}\label{eq:weighted_cycle_product}
  \prod_{X\in\gamma}\frac{a_X}{b_X}=1
\end{equation}
for every directed cycle $\gamma$ of the underlying dual graph.
\end{lemma}
\begin{proof}
If $v>0$ satisfies~\eqref{eq:weighted_flow}, then
$v_{\rho_{\rm cons}(X)}/v_{\rho_{\rm prod}(X)}=a_X/b_X$; multiplying these
identities around a directed cycle gives~\eqref{eq:weighted_cycle_product}.
Conversely, fix a reaction-node $\rho_0$ and set $v_{\rho_0}=1$.  For any
reaction-node $\rho$, choose a directed path from $\rho_0$ to $\rho$ and define
$v_\rho$ as the product of $a_X/b_X$ along that path.  This value is independent
of the chosen path: append a directed path from $\rho$ back to $\rho_0$ to two
candidate paths; each resulting closed directed walk decomposes into directed
cycles, so~\eqref{eq:weighted_cycle_product} makes the two products equal.
The resulting positive vector satisfies~\eqref{eq:weighted_flow} on every
edge, so $Mv=0$.
\end{proof}
Condition~\eqref{eq:weighted_cycle_product} has a counterpart in the
species--reaction graph of Craciun and Feinberg~\cite{CraciunFeinberg2006}.
In that graph, a directed cycle $\gamma$ of the dual corresponds to a cycle
that alternates between the reactions and substrates of $\gamma$, and its
edges carry the stoichiometric coefficients $a_X$ and $b_X$.  Such a cycle is
called an \emph{s-cycle} when the products of the coefficients on its two
alternating sets of edges agree, that is, when
$\prod_{X\in\gamma}a_X=\prod_{X\in\gamma}b_X$.  Thus
Lemma~\ref{lem:weighted_current_criterion} says that a weighted cycle has a
strictly positive current exactly when every directed cycle of its dual is an
s-cycle.  In~\cite{CraciunFeinberg2006} s-cycles enter conditions for
injectivity; here they characterize the existence of a positive current.

\begin{lemma}[Dual-cycle laws for weighted cycles]\label{lem:weighted_law}
Let $\gamma$ be a directed cycle in the dual of a weighted cycle and define
$p^\gamma\in\mathbb R^{\mathcal S}$ by
\begin{equation}\label{eq:weighted_p}
  p^\gamma_X=\frac{1}{b_X\,v_{\rho_{\rm cons}(X)}}
  \quad (X\in\mathcal S(\gamma)),
  \qquad
  p^\gamma_X=0\quad (X\notin\mathcal S(\gamma)).
\end{equation}
Then $p^\gamma\in\ker M^{\mathsf T}$, and its lift is the nonnegative
conservation law
\begin{equation}\label{eq:weighted_lifted_law}
  L_\gamma=\sum_{X\in\mathcal S(\gamma)}p^\gamma_XX
  +\sum_{\rho\in\mathcal R(\gamma)}\frac{1}{v_\rho}C_\rho .
\end{equation}
\end{lemma}
\begin{proof}
Let $\rho\in\mathcal R(\gamma)$, with incoming edge $X$ and outgoing edge $Y$
on $\gamma$.  Then $\rho=\rho_{\rm cons}(X)=\rho_{\rm prod}(Y)$, and
by~\eqref{eq:weighted_flow} applied to $v$,
\begin{equation}\label{eq:weighted_check}
  \sum_Z p^\gamma_Zt_\rho(Z)=b_Xp^\gamma_X=\frac{1}{v_\rho},
  \qquad
  \sum_Z p^\gamma_Zh_\rho(Z)=a_Yp^\gamma_Y
  =\frac{a_Y}{b_Y v_{\rho_{\rm cons}(Y)}}=\frac{1}{v_\rho}.
\end{equation}
Reactions off $\gamma$ meet no edge of $\gamma$.
Thus~\eqref{eq:composite_conservation_balance} holds, and
$\tau_\rho(p^\gamma)=1/v_\rho$ for $\rho\in\mathcal R(\gamma)$ and $0$
otherwise.
\end{proof}
For Example~\ref{ex:2A_B},~\eqref{eq:weighted_p} gives $p_A=1/2$ and
$p_B=1$, so $L_\gamma=\tfrac12A+B+C_1+C_2$.  This is half of
$A+2B+2C_1+2C_2$, the total amount of $A$ counted in all its forms.

\begin{theorem}[Weighted cycles]\label{thm:weighted}
Theorem~\ref{thm:marked_main} holds with ``cycle'' replaced by ``weighted
cycle''.
\end{theorem}
\begin{proof}
The proof of Section~\ref{sec:cycleproof} goes through with four changes.
First, the common flux~\eqref{eq:cycle_flux_identity} becomes $J=\lambda v$,
so $c_\rho^{(2)}/c_\rho^{(1)}=\lambda_2/\lambda_1$ still holds and the
comparison identities and Lemma~\ref{lem:tailpos} are unchanged.  Second, in the
uniqueness proof the dual edge $X$ receives the weight $b_X\mu_X$, which has
the sign of $\mu_X$; the incoming sum at $\rho$ is again $\tau_\rho(\mu)$.
Third, Corollary~\ref{cor:dcycle_law} is replaced by Lemma~\ref{lem:weighted_law}.
Fourth, in the proof of Proposition~\ref{prop:injective}, $\ker M=\mathbb Rv$
gives $\Delta J=\delta v$ and $\Delta c_\rho=\delta v_\rho/k_\rho$, and the
dual edge $X$ again receives the weight $b_X\mu_X$.
The remaining arguments use only that $L_\gamma$ is a nonnegative conservation
law whose support is $\mathcal S(\gamma)\cup\{C_\rho:\rho\in\mathcal
R(\gamma)\}$, which still holds.
\end{proof}

By Theorem~\ref{thm:weighted}, every HMM realization of the hypergraph of
Example~\ref{ex:2A_B} has exactly one positive steady state in each positive
compatibility class and is persistent.  The cycle condition is therefore not
necessary for the conclusion of Theorem~\ref{thm:marked_main}.

Proposition~\ref{prop:positive_classes} holds for a weighted cycle with $L_\gamma$ as in
Lemma~\ref{lem:weighted_law}: there $p\ge0$ in $\ker M^{\mathsf T}$ gives the
ordinary circulation $X\mapsto b_Xv_{\rho_{\rm cons}(X)}p_X$, whose
directed-cycle decomposition yields the weighted laws.

\subsection{Chains of intermediates and direct reactions}
The elimination of intermediate species at steady state is developed
systematically by Feliu and Wiuf~\cite{FeliuWiuf2013}; see also
Marcondes de Freitas, Feliu, and Wiuf~\cite{MarcondesDeFreitasFeliuWiuf2017}
for preservation of persistence and boundary-steady-state criteria under
removal of intermediates and catalysts.  Here the cycle laws make it possible
to obtain a uniform conclusion for every enzyme assignment.

\begin{proposition}[Longer intermediate chains and direct reactions]
\label{prop:chains}
Theorems~\ref{thm:marked_main} and~\ref{thm:weighted} remain true if each
reaction $\rho$ is realized either by a chain
\begin{equation}\label{eq:chain_module}
  E_{\sigma(\rho)}+t_\rho
  \;\rightleftarrows\; C_{\rho,1}
  \;\rightleftarrows\; C_{\rho,2}
  \;\rightleftarrows\;\cdots
  \;\rightleftarrows\; C_{\rho,\ell_\rho}
  \;\longrightarrow\; E_{\sigma(\rho)}+h_\rho,
  \qquad \ell_\rho\ge1,
\end{equation}
with positive forward rate constants and nonnegative backward rate
constants, so that the binding step and any internal step may be
irreversible, or, when $\sigma(\rho)=0$, by the single
mass-action reaction $t_\rho\to h_\rho$.
\end{proposition}
\begin{proof}
\noindent\emph{Steady-state reduction.}
Fix a chain module.  Let
$\alpha_i,\beta_i$ be the rate constants of
$C_{\rho,i}\to C_{\rho,i+1}$ and
$C_{\rho,i+1}\to C_{\rho,i}$, and let $\kappa_\rho$ be the rate constant of
the final step.  At a steady state, the net flux through every link equals
$J_\rho\coloneqq \kappa_\rho c_{\rho,\ell_\rho}$:
\begin{equation}\label{eq:chain_steady}
  k_\rho^+e_{\sigma(\rho)}x^{t_\rho}-k_\rho^-c_{\rho,1}=J_\rho,
  \qquad
  \alpha_ic_{\rho,i}-\beta_ic_{\rho,i+1}=J_\rho
  \quad(1\le i<\ell_\rho).
\end{equation}
Solving from $c_{\rho,\ell_\rho}=J_\rho/\kappa_\rho$ backward through the chain
gives
$c_{\rho,i}=\pi_{\rho,i}J_\rho$ with constants $\pi_{\rho,i}>0$
depending only on rate constants, since $c_{\rho,i}=(J_\rho+\beta_ic_{\rho,i+1})/\alpha_i$
with $\alpha_i>0$ and $\beta_i\ge0$, and then
$J_\rho=k_\rho^*e_{\sigma(\rho)}x^{t_\rho}$ with
$k_\rho^*\coloneqq k_\rho^+/(1+k_\rho^-\pi_{\rho,1})>0$.  The substrate equations
reduce to $MJ=0$.  For two positive steady states,
$c_{\rho,i}^{(2)}/c_{\rho,i}^{(1)}=\lambda_2/\lambda_1$ for every intermediate
in a cycle, and likewise in a weighted cycle after writing $J=\lambda v$.
Thus every intermediate has logarithmic change $\theta$.  The enzyme total is
\[
  T_a
  =E_a+\sum_{\rho\in\mathcal R_a}\sum_{i=1}^{\ell_\rho}C_{\rho,i},
\]
so the enzyme-comparison proof and hence the uniqueness argument are
unchanged.

\smallskip
\noindent\emph{Conservation, existence, and boundary exclusion.}
For $p\in\ker M^{\mathsf T}$, the lifted law places coefficient
$\tau_\rho(p)$ on every $C_{\rho,i}$; an internal step, reversible or not, transfers
this weight between adjacent intermediates, so the law is conserved.
Similarly, the enzyme total $T_a$ contains every $C_{\rho,i}$ with
$\rho\in\mathcal R_a$.  Summing the lifted dual-cycle laws and enzyme totals
as in Lemma~\ref{lem:marked_conservative} therefore gives a conservation law
with positive coefficient on every species.  This makes each compatibility
class compact and supplies the existence argument.

At a boundary steady state with $\lambda=0$,
$c_{\rho,i}=\pi_{\rho,i}J_\rho=0$ for every intermediate.  The positive
enzyme totals make the free enzymes positive, and~\eqref{eq:chain_steady}
then forces $x^{t_\rho}=0$.  The zero-substrate dual-cycle argument
from Proposition~\ref{prop:no_boundary} applies verbatim.  The case
$\lambda>0$ is also unchanged because every intermediate, reactant substrate,
and free enzyme is then positive.

\smallskip
\noindent\emph{Persistence.}
Lemma~\ref{lem:compact-boundary-persistence} applies because the positive
conservation law makes compatibility classes compact and the argument above
excludes boundary steady states in positive classes.

For a direct reaction $t_\rho\to h_\rho$ there is no intermediate,
$J_\rho=k_\rho x^{t_\rho}$, and $\tau_\rho(\mu)=\theta$.  The laws
$L_p$ and $L_\gamma$ carry no intermediate term for $\rho$.  At a zero-current
boundary steady state the same monomial vanishes.  These observations give
uniqueness and boundary exclusion for direct reactions; compactness,
existence, and persistence use the same substrate conservation laws.
\end{proof}

\subsection{Global convergence for weighted cycles and chains}
\begin{theorem}[Global convergence for weighted cycles]
\label{thm:weighted_convergence}
Let $\mathcal H$ be a finite weighted cycle.  Suppose that each enzyme
catalyzes at most one reaction; enzyme-free reactions are allowed.  Then, for
all positive rate constants, every trajectory that starts in a positive
stoichiometric compatibility class converges to the unique positive steady
state of its stoichiometric compatibility class.
\end{theorem}
\begin{proof}
Use the proof of Section~\ref{sec:convergence}.  In
Lemma~\ref{lem:flux_kernel},
\eqref{eq:weighted_flow} gives $\ker\Gamma=\mathbb R\hat v$ with
$\hat v_{B_\rho}=\hat v_{K_\rho}=v_\rho$.  Since $\hat v\gg0$, the kernel condition~(9)
of~\cite{AngeliDeLeenheerSontag2010} holds.  The computation in
Lemma~\ref{lem:coop} uses only that each substrate lies in one tail and one
head and that the dual graph is strongly connected, so the reaction-coordinate
system is again monotone and strongly monotone.  Persistence comes from
Theorem~\ref{thm:weighted}.  The conclusion therefore follows from
\cite[Theorem~2, case~(9)]{AngeliDeLeenheerSontag2010}, with boundary initial
conditions handled by instant positivity as before.
\end{proof}

\begin{proposition}[Intermediate chains and direct reactions]
\label{prop:chains_convergence}
The conclusions of Theorem~\ref{thm:distinct_convergence} and
Theorem~\ref{thm:weighted_convergence} remain valid when any modules are
replaced by the chains of intermediates or direct reactions in
Proposition~\ref{prop:chains}.
\end{proposition}

\begin{proof}
A direct reaction may have a substrate $A$ in both $t_\rho$ and $h_\rho$,
as in $A+B\to A+C$, which violates hypothesis~H1
of~\cite{AngeliDeLeenheerSontag2010}.  Such an $A$ occurs in no other
reaction, because it lies in one tail and one head, and flow conservation at
$A$ with $v_\rho>0$ gives equal coefficients of $A$ in $t_\rho$ and $h_\rho$.
Hence $x_A$ is constant on each compatibility class, and positive on a
positive class.  Absorbing its concentration into the rate
constant of $\rho$ and deleting $A$ gives a weighted cycle on the remaining
substrates whose dual graph differs only by the loop at $\rho$.  If $\mathcal H$
is a single reaction with $t_\rho=h_\rho$, every concentration is constant and
there is nothing to prove.  We may therefore assume that no direct reaction
has a substrate on both sides.

The graphical conditions and the kernel condition
of~\cite{AngeliDeLeenheerSontag2010} are invariant under successive removal of
intermediates~\cite[Theorems~2--3 and Corollary~1]{MarcondesDeFreitasWiufFeliu2017Convergence}.
For the present family the required check is immediate.  Splitting each chain
into elementary steps, the intermediate rows force every step of the
$\rho$-chain to have the same flux, and the substrate rows then give $MJ=0$.
Thus the positive kernel generator repeats $v_\rho$ on every step.  Each
intermediate, enzyme, and substrate occurs on opposite sides of the two rates
in which it appears, so the reaction-coordinate Jacobian is cooperative; the
arrows along the chains and the strongly connected dual make it irreducible.  An irreversible step removes only the backward arrow at that step.  The
forward arrows along each chain, from the binding step to the final step,
remain, and so do the arrows from the final step of $\rho$ to the binding
step of every $\rho'$ with a dual edge $\rho\to\rho'$; since the dual is
strongly connected, these arrows alone form a strongly connected graph.
Direct reactions are chains of one step.  Persistence follows from
Proposition~\ref{prop:chains}, so
\cite[Theorem~2, case~(9)]{AngeliDeLeenheerSontag2010} applies.
\end{proof}

\begin{corollary}[Known special cases]\label{cor:special_cases}
\begin{enumerate}
\item Let $X_1\to X_2\to\cdots\to X_m\to X_1$, $m\ge2$, be a transformation
  cycle in which each reaction has its own enzyme and is realized by a
  chain~\eqref{eq:chain_module}.  Each positive compatibility class has a
  unique positive steady state, and every trajectory in the class converges
  to it.  This is~\cite[Theorem~4.1]{EithunShiu2017} without product
  inhibition.  It contains the theorem of Rao~\cite{Rao2017}, the processive
  multisite results of~\cite{ConradiShiu2015}, and, for $m=2$ with one
  intermediate in each chain, the one-site futile cycle~\cite{AngeliSontag2008}.
\item In the cyclic double phosphorylation network
  of~\cite{SuwanmajoRameshKrishnan2020}, under any assignment of kinases and
  phosphatases, each positive compatibility class has a unique positive
  steady state and the network is persistent.  With distinct kinases and
  distinct phosphatases, every trajectory in a positive class converges to
  this steady state.
\item In the mixed mechanism~\cite{SuwanmajoKrishnan2015}
  \begin{equation}\label{eq:mixed}
  \begin{aligned}
    &S_0+K\rightleftarrows S_0K\longrightarrow S_1K\longrightarrow S_2+K,\\
    &S_2+F\rightleftarrows S_2F\longrightarrow S_1+F,\qquad
     S_1+F\rightleftarrows S_1F\longrightarrow S_0+F,
  \end{aligned}
  \end{equation}
  each positive compatibility class has a unique positive steady state and
  no boundary steady state, and the network is persistent;
  compare~\cite[Proposition~3.1]{ConradiMinchevaShiu2019}.
\end{enumerate}
\end{corollary}

\begin{proof}
In part~1 the skeleton is a cycle and no enzyme catalyzes two reactions, so
Proposition~\ref{prop:chains_convergence} applies.  In part~2 the four
phosphoforms form a transformation cycle realized by HMM modules, so
Theorems~\ref{thm:marked_main} and~\ref{thm:distinct_convergence} apply.  In
part~3 the skeleton is the transformation cycle $S_0\to S_2\to S_1\to S_0$,
the kinase reaction is a chain~\eqref{eq:chain_module} with two intermediates
and an irreversible internal step, and the two phosphatase reactions are HMM
modules, so Proposition~\ref{prop:chains} applies.
\end{proof}

\section{Discussion}\label{sec:discussion}

Substrate-cycle structure forces existence, uniqueness, and persistence
for every enzyme assignment (Theorem~\ref{thm:marked_main}), while the
enzyme-incidence pattern controls the interior dynamics
(Theorems~\ref{thm:distinct_convergence} and~\ref{thm:hopf_main}).  Both conclusions rest on the dual graph: its cycles supply conservation
laws with supports adapted to the comparison set or the boundary zero set
(Theorem~\ref{thm:marked_main}), and
its edges supply the sign pattern in reaction coordinates
(Theorems~\ref{thm:distinct_convergence} and~\ref{thm:hopf_main}).  Enzyme identities
enter the steady-state arguments only through the throughput relations and
enzyme totals, although they remain dynamically consequential away from
steady state.  The results suggest several research directions.

The first is to characterize the substrate hypergraphs for which every HMM
realization has exactly one positive steady state in each positive class.
The weighted cycles of Section~\ref{sec:extensions} show that the cycle
hypothesis is not necessary, and Example~\ref{ex:ray_counterexample} shows that
a one-dimensional current cone is not sufficient.  One route is to
replace dual cycles by a more general family of nonnegative substrate
conservation laws whose supports reproduce the sign arguments used here.

The second is to determine which enzyme-incidence patterns, on which
substrate hypergraphs, force convergence for all rate constants, and which
admit periodic or other non-equilibrium recurrent behavior.
For complex-balanced mass-action systems, persistence combines with the
Horn--Jackson Lyapunov function to give global convergence; version 3 of
Craciun's preprint~\cite{Craciun2026GAC} presents a proof of the Global Attractor
Conjecture by establishing the required persistence.  For substrate cycles,
persistence and uniqueness hold for every enzyme assignment
(Theorem~\ref{thm:marked_main}), yet convergence can fail when two reactions
share an enzyme (Theorem~\ref{thm:hopf_main}), as in the oscillations that
Suwanmajo, Ramesh, and Krishnan~\cite{SuwanmajoRameshKrishnan2020} found when
a kinase or a phosphatase is shared.  For transformation cycles, Theorems~\ref{thm:distinct_convergence}
and~\ref{thm:hopf_main} and Proposition~\ref{prop:convergent} suggest an
answer.

\begin{conjecture}\label{conj:assignments}
On a transformation cycle realized by HMM modules, with enzyme-free reactions
allowed, some positive rate constants and some positive compatibility class
admit a trajectory that does not converge to the positive steady state of
that class if and only if some enzyme catalyzes at least two reactions but
not all of them.
\end{conjecture}

The ``only if'' direction follows from Theorem~\ref{thm:distinct_convergence}
and Proposition~\ref{prop:convergent}, and
Theorem~\ref{thm:hopf_main} proves the conjecture for cycles with at most three
reactions.  For four reactions, a numerical search found an unstable positive
steady state for each of the eleven assignments, up to rotation and
relabeling of enzymes, in which some enzyme catalyzes at least two but not
all reactions (Appendix~\ref{app:four}).  The cyclic network of~\cite{SuwanmajoRameshKrishnan2020}
realizes two of them, $(E_1,E_1,E_2,E_2)$ with a common
kinase and a common phosphatase, and $(E_1,E_1,E_2,E_3)$ with two kinases and
a common phosphatase.  For its shared-kinase, shared-phosphatase
case, Conradi and
Mincheva~\cite{ConradiMincheva2024} proved Hopf bifurcations in the
subnetwork with irreversible substrate binding.  For the case
of two kinases and a common phosphatase, the oscillations
of~\cite{SuwanmajoRameshKrishnan2020} are so far only computed, and a proof
is a concrete open problem; so is the conjecture
of~\cite{ConradiMinchevaShiu2019} that the mixed mechanism~\eqref{eq:mixed}
has a stable periodic orbit whenever its steady state is unstable.

The third is to extend Theorems~\ref{thm:marked_main}
and~\ref{thm:distinct_convergence} to other catalytic mechanisms.  The
detailed mechanisms from which a substrate skeleton can be read off, built
from a reactant block and a product block, are described
in~\cite[Appendix~A.1.4]{JoshiNguyen2026Hypergraph}.  Product inhibition, in which the products rebind the
enzyme before release, is the first test case; Eithun and
Shiu~\cite{EithunShiu2017} proved global convergence for transformation cycles
with distinct enzymes and product inhibition.  Bifunctional enzymes, in
which an intermediate acts as an enzyme for another
reaction~\cite{JoshiNguyen2024,JoshiNguyen2026Hypergraph,JoshiNguyenJohnston2026},
and cascades, in which the product of one enzymatic layer acts as an enzyme in
a downstream layer~\cite{HuangFerrell1996,Kholodenko2000,QiaoEtAl2007}, are
further mechanisms to consider.

The fourth is to determine, for each dynamical property, whether it depends
only on the substrate hypergraph.  Together with~\cite{JoshiNguyen2026Hypergraph}, the present results suggest seeking
correspondences
\begin{equation}\label{eq:program}
  \text{substrate-hypergraph structure}
  \quad\Longrightarrow\quad
  \text{dynamical properties of its enzymatic realizations},
\end{equation}
for existence and multiplicity of positive steady states, boundary behavior,
concentration robustness, stability, bifurcations, and oscillation.  The
results here sort such properties into two kinds.  Existence, uniqueness, and
persistence of positive steady states depend only on the substrate
hypergraph.  Global convergence and oscillation depend in addition on the
enzyme-incidence pattern: distinct enzymes force convergence, while a shared
enzyme can produce a Hopf bifurcation.  Concentration robustness
in~\cite{JoshiNguyen2026Hypergraph} arises from the interaction of currents
with bifunctional enzymes, so it also depends on enzyme data.  For
post-translational modification systems, steady-state parametrizations
\cite{ThomsonGunawardena2009,FeliuWiuf2013PTM} provide a complementary route, and the topology--dynamics correspondences for
autocatalytic networks in~\cite{JoshiCraciun2021} a related structural
perspective.
A complementary direction is to study the prevalence of these properties in
random substrate-hypergraph models, paralleling results for multistationarity
and absolute concentration robustness in random reaction
networks~\cite{JoshiKaihnsaNguyenShiu2023}.

The fifth is a theory of inheritance and lifting at the level of substrate
hypergraphs.  The inheritance viewpoint in
reaction-network theory originated with the atoms-of-multistationarity
program of Joshi and Shiu for fully open networks~\cite{JoshiShiu2013}, which identified
conditions under which multistationarity of an embedded network lifts to a
larger network.  Banaji and coauthors showed that, under suitable hypotheses, a network's capacity for
nondegenerate multistationarity or periodic orbits is preserved by several enlargement
operations~\cite{BanajiPantea2018,Banaji2018}, including rank-preserving
additions of species and reaction
splittings~\cite{BanajiBorosHofbauer2022,Banaji2023}; related results
preserve local bifurcations under transversality
assumptions~\cite{BanajiBorosHofbauer2025}.  Conradi and
Shiu~\cite{ConradiShiu2018} asked whether Hopf bifurcations are preserved when
reactions or intermediates are added, and Conradi and
Mincheva~\cite{ConradiMincheva2024} used Banaji's result on the inheritance of oscillation~\cite{Banaji2018} to show that a
stable limit cycle of the subnetwork with irreversible binding, if one exists,
yields a stable limit cycle of the full network.
Theorem~\ref{thm:marked_main} and Proposition~\ref{prop:chains} provide first
instances of such lifting for substrate hypergraphs: existence, uniqueness, and persistence pass from
the substrate cycle to every enzyme assignment and every intermediate-chain
length, and Corollary~\ref{cor:special_cases} shows that the known
uniqueness and convergence results for processive, cyclic, and mixed
phosphorylation mechanisms are instances of it.  For systems without inflows or outflows, a natural next step is to seek cyclic decompositions
of suitably strongly connected substrate hypergraphs and to determine which
cycle-level conclusions survive when the pieces are recombined.  Such a theory
could yield inheritance and steady-state-counting results and, because enzyme
sharing already produces periodic orbits on a single cycle, criteria for when
limit cycles persist under composition.

The sixth is to treat open systems.  Inflows and outflows
destroy the substrate conservation laws that underlie the compactness,
boundary-exclusion, and persistence arguments used here, so open substrate
hypergraphs will require different tools.

The present results identify the substrate cycle as a basic unit of
inheritance.  The next test is whether composition, opening, and catalytic
refinement preserve, or systematically alter, the properties found here.  A
successful theory would make substrate hypergraphs a predictive bridge from
coarse biochemical structure to the dynamics of entire families of enzymatic
systems.

\section*{Acknowledgments}
B.~Joshi was supported by an AMS--Simons Research Enhancement Grant for
Primarily Undergraduate Institution Faculty.

\section*{Acknowledgment of AI use}
The authors used Claude (Anthropic) and ChatGPT with Codex (OpenAI) to review
the manuscript for mathematical correctness and clarity, suggest and implement
revisions to the text and some proofs, locate and check references, and verify
symbolic and numerical computations.  Claude was also used to write the scripts for the
numerical work: the exact checks, the interval-arithmetic certification of the
first Lyapunov coefficient, and the simulations of
Section~\ref{sec:instability} and Appendix~\ref{app:dynamics_computation},
including Figure~\ref{fig:hopf_cycle}; and the four-reaction search in
Appendix~\ref{app:four}.  The authors checked all AI-assisted
work and take full responsibility for the paper's content.

\section*{Code availability}
The scripts that reproduce the computations of Section~\ref{sec:instability},
Figure~\ref{fig:hopf_cycle}, and Appendices~\ref{app:dynamics_computation}
and~\ref{app:four}, together with their requirements files, are
archived on Zenodo, doi:10.5281/zenodo.23777685.

\appendix

\section{Two lemmas on directed graphs}\label{app:graph}

\begin{lemma}[Backward walk]\label{lem:backward_walk}
Let $G$ be a finite directed graph and $F$ a set of its edges such that every
vertex is the head of at least one edge of $F$.  Then $F$ contains a directed
cycle.
\end{lemma}
\begin{proof}
Start at any vertex $v_0$.  Given $v_i$, choose an edge of $F$ with head $v_i$
and let $v_{i+1}$ be its tail.  Let $j$ be the first index for which
$v_j\in\{v_0,\ldots,v_{j-1}\}$, say $v_j=v_i$.  The edges chosen at steps
$i,\ldots,j-1$, traversed in reverse order, form a directed cycle in $F$.
\end{proof}

\begin{lemma}[A nonnegative directed cycle]
\label{lem:nonnegative_edge_cycle}
Let $G$ be a finite strongly connected directed graph.  Give each edge $\alpha$
a real weight $\eta_\alpha$, not all zero, and suppose that at every vertex the
sum of the weights of the incoming edges is nonnegative.  Then some directed
cycle has all edge weights nonnegative and at least one edge weight positive.
\end{lemma}
\begin{proof}
Suppose, for contradiction, that no directed cycle has all edge weights
nonnegative and at least one edge weight positive.  Some edge has positive
weight: otherwise every weight is nonpositive, every incoming sum is then
zero, and every weight vanishes.  Let $G_{\ge0}$ be the spanning subgraph of
edges of nonnegative weight.  A positive edge with both endpoints in one
strongly connected component of $G_{\ge0}$ would lie on a directed cycle of
$G_{\ge0}$, which is excluded.  Hence every positive edge joins two different
strongly connected components of $G_{\ge0}$, and $G_{\ge0}$ has at least two
components.  The edges of $G_{\ge0}$ between distinct components form an
acyclic directed graph on the components, so some component $Q$ has no edge of
$G_{\ge0}$ entering it from another component.  Since $G$ is strongly
connected and $Q$ is not all of $G$, at least one edge of $G$ enters $Q$ from
outside, and every such edge has negative weight.  Every edge with both
endpoints in $Q$ has nonpositive weight, since every positive edge was shown
above to join two distinct strongly connected components.  Summing the
incoming weights over the vertices of $Q$
therefore gives a negative number, which contradicts the hypothesis.
\end{proof}

\section{A fixed-point lemma}\label{app:fixedpoint}
\begin{proof}[Proof of Lemma~\ref{lem:compact_equilibrium}]
Let $P_K$ denote Euclidean projection onto the closed convex set $K$, and fix
$\varepsilon>0$.  The continuous map
\[
  F(u)\coloneqq P_K\bigl(u+\varepsilon f(u)\bigr),\qquad u\in K,
\]
maps $K$ into itself, so Brouwer's fixed-point theorem gives $u^*\in K$ with
$F(u^*)=u^*$.  The defining property of Euclidean projection onto a convex set
then gives
\[
  \langle f(u^*),\tilde u-u^*\rangle\le0
  \qquad(\tilde u\in K).
\]
Let $u(t)$ be the solution starting at $u^*$.  Positive invariance gives
$u(t)\in K$ for $t\ge0$; because the trajectory remains in the compact set
$K$, the local solution extends forward for all time.  Substitute $\tilde u=u(t)$,
divide by $t>0$, and let $t\downarrow0$.  Since $(u(t)-u^*)/t\to f(u^*)$, we
obtain $\|f(u^*)\|^2\le0$, hence $f(u^*)=0$.
\end{proof}

\section{Instability and Hopf computations}\label{app:dynamics_computation}

For the family~\eqref{eq:hopf-family}, the characteristic polynomial in the
proof of Proposition~\ref{prop:hopf} has coefficients
\begin{equation}\label{eq:hopf-coefficients}
\begin{aligned}
 a_1&=\tfrac{1}{50}(110s+122851),&
 a_2&=\tfrac{1}{250}(1346561s+72272170),\\
 a_3&=\tfrac{48}{125}(1533793s+8905715),&
 a_4&=\tfrac{3}{10}(7851979s+39361895),\\
 a_5&=63641880(s+5).
\end{aligned}
\end{equation}
Its Hurwitz determinants are $\Delta_1=a_1$,
\begin{align}
 \Delta_2&=\tfrac{11}{12500}\bigl(13465610s^2+15092190701s+803269174970\bigr),
 \label{eq:hurwitz2}\\
 \Delta_3&=\tfrac{27}{3125000}\bigl(806460136003940s^3
   +907107963131768904s^2\notag\\
 &\qquad+51801090005266991525s+271628221015418837525\bigr),
 \label{eq:hurwitz3}\\
 \Delta_4&=\tfrac{9}{31250000}\,g(s),\qquad \Delta_5=a_5\Delta_4,
 \label{eq:hurwitz4}
\end{align}
where
\begin{equation}\label{eq:hopf-q}
\begin{aligned}
 g(s)={}&42886676966934711283340s^4+47754860538759245982248644s^3\\
 &+2210683395785214746557065495s^2-16047002691246038424525453350s\\
 &-129527231666398742587344911125.
\end{aligned}
\end{equation}
For the enzyme-free family~\eqref{eq:hopf-free-family} of
Proposition~\ref{prop:hopf-free}, in the same coordinates
$(x_1,x_2,c_1,c_2,c_3)$, the coefficients are
\begin{equation}\label{eq:hopf-free-coefficients}
\begin{aligned}
 a_1&=\tfrac{1}{50}(100s+122801),&
 a_2&=\tfrac{1}{50}(244852s+14332133),\\
 a_3&=\tfrac{3}{50}(8942417s+52627285),&
 a_4&=\tfrac{3}{10}(7535371s+37778855),\\
 a_5&=63635760(s+5),
\end{aligned}
\end{equation}
and the Hurwitz determinants are $\Delta_1=a_1$,
\begin{align}
 \Delta_2&=\tfrac{1}{2500}\bigl(24485200s^2+30159921202s+1752106171783\bigr),
 \notag\\
 \Delta_3&=\tfrac{27}{62500}\bigl(12143338898800s^3
   +15003846923177013s^2\notag\\
 &\qquad+927168218480716457s+4966256087694059160\bigr),
 \notag\\
 \Delta_4&=\tfrac{27}{625000}\,\tilde g(s),\qquad \Delta_5=a_5\Delta_4,
 \notag
\end{align}
where
\begin{equation}\label{eq:hopf-free-q}
\begin{aligned}
 \tilde g(s)={}&203863264999818604400s^4+249048449576344929284869s^3\\
 &+12063041773993324064613186s^2-129239390345935782322640315s\\
 &-916553932488285472722472600.
\end{aligned}
\end{equation}

The accompanying script \texttt{hypcycle\_dynamics.py} reproduces the
calculations in Section~\ref{sec:instability} without importing or executing
any other source file.  Its \texttt{exact} command uses rational SymPy
arithmetic to verify the stated steady state and conservation totals, compute
the reduced Jacobian and characteristic polynomial, derive the Hurwitz
determinants and the polynomial $g$, isolate its unique positive root, and
assert that the instability and transversality checks pass.  The
\texttt{certify} command uses the Arb ball arithmetic of
\texttt{python-flint}~\cite{Johansson2017}.  It isolates $s^*$ in a rational
interval of width less than $10^{-60}$, evaluates the Jacobian, the second
derivatives, and $\omega^2=(a_1a_4-a_5)/(a_1a_2-a_3)$ on that interval, and
computes enclosures of the right and left eigenvectors $q$ and $p$ by fixing
their $x_2$-components and solving four of the five equations.  At $s=s^*$,
Liu's criterion gives one-dimensional kernels for $A-i\omega I$ and
$A^{\mathsf T}+i\omega I$.  The interval solver certifies that the selected
$4\times4$ minors are nonsingular; hence the $x_2$-components of the kernel
vectors are nonzero, fixing them to $1$ uniquely determines the true
eigenvectors, and the interval solves enclose them.  The enclosures are
then rescaled so that $\langle q,q\rangle=1$ and $\langle p,q\rangle=1$, the
normalization used in~\eqref{eq:lyapunov}.  The omitted equations
then follow from singularity of the two matrices; their interval residuals
are checked only as consistency tests.  The command next encloses
$h_{11}=A^{-1}B(q,\bar q)$,
$h_{20}=(2i\omega I-A)^{-1}B(q,q)$, and $\ell_1$
in~\eqref{eq:lyapunov}, and asserts $\ell_1<0$.  The
\texttt{lyapunov} command evaluates \eqref{eq:lyapunov} in floating point at high precision as an independent check.  The \texttt{simulate} command
integrates the reduced five-dimensional system at $s=1/2$ and estimates the
late-time period and $x_2$ range.  The symbolic identities and
interval-arithmetic enclosures provide rigorous computer-assisted
verification.  The floating-point eigenvalue and Lyapunov-coefficient
calculations, and the trajectory integration, provide independent numerical
checks.  A separate script,
\texttt{make\_hopf\_figure.py}, integrates the full eight-dimensional system
with the same rate constants and produces Figure~\ref{fig:hopf_cycle}.

After installing the versions in \texttt{requirements-dynamics.txt}, all
checks may be run with
\begin{center}
\texttt{python hypcycle\_dynamics.py all}.
\end{center}
The individual commands \texttt{exact}, \texttt{certify}, \texttt{lyapunov}, and
\texttt{simulate} can be used separately.  The expected leading numerical
outputs are $s^*=10.4673375245\ldots$,
$\omega=10.27548989\ldots$, a first Lyapunov coefficient $\ell_1=-2.68605139416826643\ldots\times10^{-3}$, and, for $s=1/2$, a period near $1.28$ with $x_2$
ranging approximately from $19$ to $66$.  The option
\texttt{--variant enzyme-free}, placed before the command, runs the same
checks for Proposition~\ref{prop:hopf-free}; the expected outputs are
$\tilde s^*=13.0160798148\ldots$, $\tilde\omega=10.78715319\ldots$,
$\ell_1=-2.94698791053829051\ldots\times10^{-3}$, and, for $s=1/2$, a period
near $1.33$ with $x_2$ ranging approximately from $16$ to $66$.

\section{Four-reaction transformation cycles}\label{app:four}

For four reactions there are eighteen assignments up to rotation and
relabeling of enzymes, where $0$ marks an enzyme-free reaction.  In seven of
them no enzyme catalyzes two or more reactions, or one enzyme catalyzes all
four:
\begin{itemize}\raggedright
\item each enzyme catalyzes at most one reaction: $(0,0,0,0)$, $(E_1,0,0,0)$,
  $(E_1,E_2,0,0)$, $(E_1,0,E_2,0)$, $(E_1,E_2,E_3,0)$, $(E_1,E_2,E_3,E_4)$;
\item one enzyme catalyzes all four reactions: $(E_1,E_1,E_1,E_1)$.
\end{itemize}
These cases are covered by Theorem~\ref{thm:distinct_convergence} and
Proposition~\ref{prop:convergent}, respectively: for all positive rate
constants, every trajectory in a positive compatibility class converges to
its positive steady state.  In the remaining eleven, some enzyme catalyzes at
least two but not all reactions:
\begin{itemize}\raggedright
\item two enzyme-free reactions: $(E_1,E_1,0,0)$, $(E_1,0,E_1,0)$;
\item one enzyme-free reaction: $(E_1,E_1,E_1,0)$, $(E_1,E_1,E_2,0)$,
  $(E_1,E_1,0,E_2)$, $(E_1,E_2,E_1,0)$;
\item no enzyme-free reaction: $(E_1,E_1,E_1,E_2)$, $(E_1,E_1,E_2,E_2)$,
  $(E_1,E_1,E_2,E_3)$, $(E_1,E_2,E_1,E_2)$, $(E_1,E_2,E_1,E_3)$.
\end{itemize}
These cases were tested numerically.  For each assignment, $10^4$ samples
were drawn, each consisting of a positive steady state and unbinding and
catalytic rate constants chosen log-uniformly; the binding rate constants
were then determined by the steady-state equations, and the eigenvalues of
the Jacobian on the compatibility class were computed (script
\texttt{hypcycle\_conj4\_full.py}).  For each of the eleven assignments, some
sample gave an unstable positive steady state, with an eigenvalue of positive
real part.

\end{document}